# ARE Method[1]: Orbital Decompositions and Dihedral Cancellations for Determinants

Ramón Moya

School of Mathematics, Autonomous University of Santo Domingo (UASD), Dominican Republic

rmoya07@uasd.edu.do

**ORCID:** 0009-0001-1601-4699

**Document scope.**

This manuscript is a technical reference document for the ARE Method. It is intentionally expansive and includes definitions, proofs, examples, visualization methods, implementation notes, computational tests, and supplementary explanations. Its purpose is to present the full development of the framework in a form suitable for future study, verification, reuse, and extension.

A shorter article version may be extracted from this document. However, the present manuscript should not be read as a compact journal article, but as an extended research monograph documenting the complete mathematical framework.

**Abstract**

We introduce the **ARE Method (Action, Rectification, and Structure)**, a research framework for organizing the Leibniz expansion of the determinant before summation. The method partitions the $n!$ Leibniz terms into cyclic orbits under the right action of the cyclic group $C_n$ on $S_n$. Each orbit has size n and admits a canonical representative satisfying $\sigma^\star(1) = 1$.

The first level of the framework is the cyclic action, which decomposes the full symmetric group into $(n-1)!$ orbital families. The second level is canonical rectification: after a suitable column permutation, the monomials of a given orbit appear as parallel diagonal families with consecutive offsets. The third level is dihedral pairing, which associates each orbit with a companion orbit and provides a precise setting for studying cancellations under additional weight-invariance hypotheses.

[1]The acronym ARE stands for Action, Rectification, and Structure,the three pillars of the method. The author dedicates this work to his children, Arlette and Ramón Eduardo, whose initials share the same acronym.

The framework does not replace Gaussian elimination or LU decomposition and does not improve the general asymptotic complexity of determinant computation. Its contribution is structural, conceptual, and pedagogical: it explains the cyclic organization behind Sarrus-type diagrams, clarifies why fixed-width Sarrus extensions fail for $n \geq 4$, and provides a systematic language for orbital signs, canonical rectification, and conditional dihedral cancellations.

The work also connects the cyclic organization of Leibniz terms with the Discrete Fourier Transform: cyclic rotation corresponds to circular translation, and the determinant appears as the zero-frequency component of an associated orbital distribution. We provide formal definitions, proofs, examples in small dimensions, several visualization methods, and reproducible computational verification.

**Reproducibility.** Code, data, and scripts for regenerating figures/tables are available from Zenodo (DOI: 10.5281/zenodo.17423738). The repository includes sarrus_gui_v12.py (interactive GUI), orbital_methods.py (pure mathematical functions), and test_matrices.py (verification of 147 matrices).

This approach formalizes and extends Arschon's pioneering geometric ideas [1] by means of the algebraic framework of orbits under the action of the cyclic group $C_n$ on $S_n$.

# 1 Introduction

## 1.1 Motivation and Context

The determinant of an $n \times n$ matrix A is defined by the Leibniz expansion as a signed sum over all permutations $\sigma \in S_n$. Each term is a product of n entries, one from each row and one from each column. For $n = 3$, Sarrus's rule organizes the six Leibniz terms into three positive and three negative visual diagonals. For $n \geq 4$, however, a direct fixed-width extension of this visual scheme cannot represent all n! Leibniz terms as a single family of Sarrus-type diagonals.

The correct question is not merely how to extend the Sarrus diagram, but which algebraic structure underlies it and how that structure can be generalized to arbitrary dimension. The answer we present is the **ARE Method**: the $n!$ Leibniz terms are not a flat set, they are a union of $(n-1)!$ orbits of size $n$ under the action of counterclockwise rotation $\sigma \circ \rho^r$ of the cyclic group $C_n$.

The counterclockwise action has a crucial property: it is algebraically identical to the cyclic translation that generates the rows of a circulant matrix and underlies the Discrete Fourier Transform. Therefore, the determinant is recoverable as the zero-frequency component. $G_0$ of a complex determinant vector $D(A) \in \mathbb{C}^n$, where each mode $G_k = \sum_\sigma m_\sigma(A) \cdot \omega^{k\varphi(\sigma)}$ is a Fourier coefficient of the monomial distribution. The cyclic action used here is formally compatible with the circular translations that underlie circulant matrices and the Discrete Fourier Transform. This does not mean that the determinant of a general matrix is computed by a DFT. Rather, it means that once the Leibniz terms are grouped into cyclic orbital classes, one may form orbital distributions whose Fourier coefficients encode modal information. In that sense, the ordinary determinant corresponds to the zero-frequency component of the associated orbital distribution.

This article establishes the foundations of the ARE Method and demonstrates that the base monomials, all with $\sigma(1) = 1$, which canonically anchors the representative of each orbit, generate $S_n$ complete by counterclockwise rotation and canonical rectification. The supplementary articles develop: (i) the determinant vector $D(A)$ as a spectral invariant [13]; (ii) the TDet orbital hyperdeterminants for tensors of order 3 (a companion paper in preparation).

### 1.2 Our Contribution

**Impossibility of a single fixed-width Sarrus-type rectification**

Sarrus's rule for $n = 3$ uses an extended matrix of width 5, obtained by repeating the first two columns of the original matrix, in order to display the six Leibniz terms as three descending and three ascending diagonal products. A natural question is whether the same visual principle can be extended to arbitrary n by means of a single fixed-width diagram.

The following result clarifies the precise limitation. It does not say that determinants cannot be visualized for $n \geq 4$. It says that a single Sarrus-type rectification, based on one fixed column arrangement or one orbital-pair display, cannot simultaneously represent all n! Leibniz terms as parallel diagonal families.

**Theorem 1.1 (Limitation of a single fixed-width Sarrus-type rectification).**

Let A be an $n \times n$ matrix. For $n \geq 4$, no Sarrus-type scheme based on a single fixed-width representation, or equivalently on a single orbital-pair rectification, can display all $n!$ Leibniz terms as parallel diagonal families arising from one column arrangement.

Proof.

A canonical rectification associated with a cyclic orbit $\Omega(\sigma)$ aligns the $n$ monomials $\Omega(\sigma)$ as diagonal families with consecutive offsets. Under the corresponding dihedral pairing, the same rectification also displays the $n$ monomials of its companion orbit as antidiagonal families. Thus, one rectification displays one orbital companion pair, namely 2n Leibniz terms.

For $n = 3$, the symmetric group $S_3$ has two cyclic orbits under the action of $C_3$. These two orbits are companions of each other. Hence a single rectification displays all $2 \cdot 3 = 6$ Leibniz terms. This is the algebraic mechanism behind the classical Sarrus rule.

For $n \geq 4$, the number of cyclic orbits is

$(n!)/n = (n-1)!$.

These orbits are grouped into companion pairs, so the number of orbital companion pairs is

$\big((n-1)!\big)/2$.

Since $\big((n-1)!\big)/2 > 1$ for $n \geq 4$, there is more than one companion pair. A single column arrangement can canonically rectify one such pair, but distinct companion pairs generally require distinct rectifications. Therefore, one fixed-width Sarrus-type diagram cannot simultaneously rectify all companion pairs and hence cannot display all $n!$ Leibniz terms as one family of parallel diagonal or antidiagonal products.

This proves the limitation of any single fixed-width Sarrus-type rectification for $n \geq 4$.

Corollary 1.2 (Consequences for Sarrus-type extensions).

(a) For $n = 3$, Sarrus's rule works because the two cyclic orbits of $S_3$ are companions. A single rectification displays both orbits and therefore covers all six Leibniz terms.

(b) For $n \geq 4$, the cyclic decomposition contains more than one companion pair. Therefore, a single rectification cannot display all Leibniz terms as parallel diagonal families.

(c) To obtain a complete visualization for $n \geq 4$ within the ARE framework, one must use separate orbital or orbital-pair rectifications. Taken together, these rectifications cover all n! Leibniz terms without omission or repetition.

(d) This limitation is not a limitation of the determinant itself, nor of visualization in general. It is a limitation of the specific Sarrus-type requirement that all terms be displayed through one fixed-width arrangement.

**Observation 1.3 (Historical and recent attempts).**

Several works have attempted to extend Sarrus's rule beyond the $3 \times 3$ case by using extended matrices, path schemes, or families of column permutations. The work closest in spirit to the present approach is that of Lorenz and Wirths, who introduce the dihedrant as a sum over a dihedral group rather than over the full symmetric group $S_n$. Their work shows that such Sarrus-type constructions naturally lead to dihedral substructures, but that the resulting dihedrant is not, in general, equal to the determinant for $n \geq 4$.

The present monograph differs from those approaches in three respects. First, it organizes the full symmetric group $S_n$, not only a dihedral subgroup, into cyclic orbits under the action of $C_n$. Second, it formulates canonical rectification as an algebraic theorem rather than only as a graphical procedure. Third, it separates the always-valid geometric pairing of orbits from the conditional algebraic cancellation of their Leibniz contributions.

Thus, Theorem 1.1 should be read narrowly: it explains why one fixed-width Sarrus-type rectification cannot cover all n! Leibniz terms for $n \geq 4$. It does not deny the value of previous visual methods, nor does it claim that all Sarrus-inspired constructions are invalid. Rather, it identifies the precise structural obstruction that motivates the orbital-pair rectifications used in the ARE Method.

**Resolution by orbital structure**

Theorem 1.1 does not state that it is impossible to visualize determinants of order $n \geq 4$. It states only that the strategy of using one fixed-width Sarrus-type rectification is insufficient.

The ARE Method resolves this by abandoning the requirement of a single global rectification. Instead, it proceeds orbitally:

(1) work with the original matrix A as the algebraic object of computation;

(2) decompose the n! Leibniz terms into $(n-1)!$ cyclic orbits under the right action of $C_n$ on $S_n$;

(3) choose one canonical representative $\sigma^\star$ for each orbit, usually the unique representative satisfying $\sigma^\star(1) = 1$;

(4) rectify each orbit separately, so that its $n$ monomials become parallel diagonal families with consecutive offsets;

(5) pair each orbit with its dihedral companion when studying dual visualization and possible cancellations.

In this sense, Sarrus's rule is not an isolated trick. It is the $n = 3$ case in which the full orbital structure collapses into a single companion pair. For $n \geq 4$, the same principle persists, but the complete determinant requires several orbital rectifications rather than one global diagram.

**Main contributions of the framework**

The ARE Method provides a structural reorganization of the Leibniz expansion. Its main contributions are the following.

Cyclic orbital decomposition. The n! Leibniz terms are partitioned into $(n-1)!$ cyclic orbits of size n under the right action of $C_n$ on $S_n$.

Canonical representatives. Each orbit contains a unique representative satisfying $\sigma(1) = 1$. This representative serves as the generator orbital for the corresponding family.

Canonical rectification. For each orbit, a suitable column permutation aligns the n orbital monomials into parallel diagonal families with consecutive offsets.

Orbital sign rule. Since $sgn(\rho) = (-1)^{(n-1)}$, the signs within an orbit are constant when n is odd and alternate when n is even.

Dihedral pairing. Each orbit admits a companion orbit under a reversal operation. This pairing is geometrically natural and explains the dual diagonal-antidiagonal structure behind Sarrus-type diagrams.

Conditional cancellation. Dihedral cancellation is not automatic. It occurs only under explicit weight-invariance hypotheses such as $W_A\big(\Phi(\tau)\big) = W_A(\tau)$, or under special structural conditions that imply this equality.

Computational reproducibility. The monograph includes algorithms, examples, and computational verification against standard determinant routines. These computations support the implementation but do not replace the mathematical proofs.

**Complexity and scope**

The ARE Method reorganizes the Leibniz expansion; it does not replace it by a polynomial-time determinant algorithm. In the general case, the number of Leibniz terms remains $n!$, and the orbital decomposition contains $(n-1)!$ orbits of size n.

Therefore, the general computational complexity remains factorial when all terms are explicitly enumerated. The method is not intended to compete with Gaussian elimination, LU decomposition, or other $O(n^3)$-type algorithms for numerical determinant computation.

Its value is structural, conceptual, visual, and pedagogical. It explains how the Leibniz terms are organized by cyclic symmetry, how Sarrus's rule fits into a general group-theoretic pattern, and how additional matrix structure may produce reductions or cancellations in special cases. Such reductions must always be stated as conditional, not as general complexity improvements.

**Five Visualization Methods Associated with the ARE Framework**

This work introduces five methods of geometric visualization for determinants $n \times n$, The first four visualizations represent the same Leibniz expansion under the main cyclic-orbital framework, although they differ in layout and purpose. The fifth method is a pedagogical additive variant; it is not the same orbit structure as the main $C_n$ action and must be distinguished notation-wise.

Although several recent works have attempted to extend Sarrus's rule [9], [11], none provide a complete theoretical framework based on group theory that unifies the geometric extensions with the underlying algebraic structure. Our approach is distinguished by establishing a rigorous correspondence between the action of the cyclic group $C_n$ and the visual representations using polylines and parallel lines.

In particular, Lorenz & Wirths [11] show that the "False Sarrus Rule" calculates the dihedrant, a quantity distinct from the determinant for, $n > 3$ and identify that it coincides with $\det(A)$ only under specific structural conditions ( $n \equiv 2,3 \pmod 4$ for antidiagonal triangular matrices). Its scheme for $4 \times 4$ use $\frac{(n-1)!}{2} = 3$ column permutations and visually anticipates our orbital construction, although without the algebraic foundation of Theorem 3.2 or the generalization to $n$ arbitrary.

**Method 1:** Standard Polylines (Appendix A.2)

**Method 2:** Rectified parallel lines (Theorem 3.2)

**Method 3:** Dihedral base-companion pairing (Theorem 4.1)

**Method 4:** Complete visualization by column rectification (successive application of Theorem 3.2)

Arschon's work [1] represents an early milestone in attempts to systematize the calculation of determinants for matrices beyond $3 \times 3$, anticipating modern computational approaches by decades. Although his geometric constructions lacked the group-theoretic formalization presented here, they demonstrated a remarkable intuition about the underlying structure. Our method of total lines provides the algebraic foundation that formalizes and generalizes these early geometric intuitions through the orbit theory of the cyclic group.

**Method 5:** Polylines by modular increment (Section 6.2)

Important: Generates a different orbit $\Omega^a(\sigma)$ under the action of $\mathbb{Z}_n$ (not equivalent to) $\Omega(\sigma)$ under $C_n$). Pedagogical value: demonstrates flexibility of the polyline method and arithmetic simplicity. Requires only $(n-1)!$ total graphs. Didactic advantage: construction without permutation compositions (accessible without advanced group theory) [2].

## 2 Preliminaries and notation

### 2.1 Groups and permutations

[2] The star block method and other alternative partitioning methods of $S_n$ They will be developed in complementary works.

Denote by $S_n$ the symmetric group of permutations of $[n] = \{1, \dots, n\}$. Let $\rho = (1\ 2 \dots n)$ the n-cycle (with the convention $\rho(i) = i + 1 \bmod n$, counterclockwise rotation) and $C_n = \langle \rho \rangle$ the cyclic group.

Let $J$ the inversion of columns (central reflection) and $D_n = \langle \rho, J \rangle$ the dihedral group. For $\sigma \in S_n$, its cyclic orbit is $\Omega(\sigma) = \{\sigma \circ \rho^r : 0 \le r < n\}$. The sign is $\operatorname{sgn}(\sigma) = (-1)^{\operatorname{inv}(\sigma)}$.

### 2.2 Offsets and families of parallel lines

In a matrix $A = (a_{i,j})$, we define $\Delta_+(i,j) \equiv (j - i) mod n$ (slope $+1$) and $\Delta_-(i,j) \equiv (i + j) mod n$ (slope $-1$). Each monomial induces a set of pairs $\{(i, \sigma(i))\}_{i=1}^{n}$, which align in families of parallel lines when rectifying the orbit.

### 2.3 Matrix Symmetry Classes

Let $A \in K^{n \times n}$ (with $K$ a field). We consider:

1. **Toeplitz matrices:** $a_{i,j} = t_{j-i}$ (constants along diagonals parallel to the main diagonal)
2. **Persymmetric matrices:** $a_{i,j} = a_{n+1-j,n+1-i}$ (symmetric with respect to the antidiagonal)
3. **Centrosymmetric matrices:** $a_{i,j} = a_{n+1-i,n+1-j}$ (symmetrical with respect to the center)

**Note:** Every centrosymmetric matrix is persymmetric, but not vice versa.

### 2.4 Leibniz's formula

The determinant is defined as:

$$det(A) = \sum_{\sigma \in S_n} \operatorname{sgn}(\sigma) \prod_{i=1}^{n} a_{i,\sigma(i)}$$

**Weight of a permutation:** We denote $W_A(\sigma) \coloneqq \prod_{i=1}^{n} a_{i,\sigma(i)}$.

## 3 Orbital decomposition and canonical rectification

### 3.1 Partitioning $S_n$ into orbits

**Theorem 3.1 (Partition into $(n-1)!$ orbits).**

The action of $C_n$ The $S_n$ composition on the right partitions $S_n$ into exactly $(n-1)!$ orbits, each of size $n$.

**Proof.**

**(1) Size of the orbits.**

For [ $\sigma \in S_n$], its stabilizer is: $\text{Stab}(\sigma) = g \in C_n : \sigma \circ g = \sigma$

We affirm that $|\text{Stab}(\sigma)| = 1$ for $n > 1$:

- If $\sigma \circ \rho^r =$with $0 \leq r < n$, then $\sigma(i) = \sigma(i + r \ mod \ n)$ for all $i \in [n]$.
- Due to injectivity σ: $if \ \sigma(i) = \sigma(j)$ then $i = j$.
- This force $i = i + r \ mod \ n$ for all $i$, which implies $r \equiv 0 \ mod \ n$.
- As $0 \leq r < n$ we conclude $r = 0$.

By the orbit-stabilizer theorem:

$$|\Omega(\sigma)| = \frac{|C_n|}{|\text{Stab}(\sigma)|} = \frac{n}{1} = n$$

**(2) Number of orbits.**

How $S_n$ partitions into disjoint orbits of size

$$n\text{:\# orbits} = \frac{|S_n|}{n} = \frac{n!}{n} = (n-1)!$$

**Corollary 3.1.** Each Leibniz term belongs to exactly one orbit, grouping the $n!$ terms into $(n-1)!$ families of $n$ terms related by rotation.

**Corollary 3.2. (Cardinality of Rectifiers).**

For each orbit $\mathcal{O}(\sigma)$, there are exactly $n$ different rectifiers, given by the family $\{\tau_\kappa \circ \sigma^{-1} : \kappa \in \{0,1, \ldots, n-1\}\}$, where $\tau_\kappa$ denotes the cyclical displacement of columns in $\kappa$ positions. All these rectifiers produce the same geometric pattern ( $n$ consecutive parallel lines), differing only in the starting diagonal. The canonical rectifier $B_{\text{can}} = \sigma^{-1}$ corresponds to $\kappa = 0$ and it is the only one that aligns the base monomial on the main diagonal.

**Definition (Canonical representative and generator orbital).**

For each cyclic orbit $\Omega(\sigma)$, we fix the unique representative $\sigma^\star$ satisfying

$\sigma^\star(1) = 1$.

This representative is called the canonical representative, or generator orbital, of the orbit. The remaining elements

$\sigma^{\star} \circ \rho^{r}, 1 \leq r < n,$

are called **orbital descendants**. Algebraically, all elements of the orbit have the same status, but the choice of $\sigma^{\star}$ provides a fixed anchor for rectification, notation, and computation. The canonical rectifier is the unique $B_{\text{can}} = (\sigma^{\star})^{-1}$, which aligns $\sigma^{\star}$on the main diagonal with consecutive offsets $0,1, \dots, n-1$.

**Example (n=4).**

Let $\sigma = [2,4,1,3]$. Its orbit is

$\{\sigma, \sigma \circ \rho, \sigma \circ \rho^{2}, \sigma \circ \rho^{3}\}$

- Based on $\sigma$:$B_{\text{can}} = \sigma^{-1}$ and the offsets (slope $+1$) are 0,1,2,3.
- Based on $\tau = \sigma \circ \rho^{2}$:$B_{\text{can}}(\tau) = \rho^{-2} \circ \sigma^{-1}$ and the offsets become 2,3,0,1; the diagram is the same, rotated.

### 3.2 Canonical rectification and pattern of signs

**Convention.** Throughout this section, expressions of the form $i + k (mod\ n)$ are interpreted in $\mathbb{Z}/n\mathbb{Z}$, with representatives in $\{1, \dots, n\}$. In particular, $n + 1 \equiv 1 (mod\ n)$.

#### 3.2.1 Canonical Rectification

**Theorem 3.2 (Canonical Rectification of a Block).** For each orbit, $\Omega(\sigma^{\star})$ there exists a canonical rectification of the matrix A that aligns the monomials of A $\Omega(\sigma)$ into families of n parallel lines with consecutive displacements. This rectification is unique up to rotation.

**Proof.**

We work with slope $+1$ (diagonal direction); the slope case $-1$ continues by reflection $A \ \mapsto \ AJ$.

**(1) Rectifier Construction**

Canonical rectification consists of multiplying by the permutation matrix P(σ) on the right:

$A^{*} \coloneqq A \cdot P(\sigma^{\star})$

where $P(\sigma^{\star})$ is the permutation matrix such that $\left(AP(\sigma^{\star})\right)\{ij\} = a\{i, \sigma^{\star}(j)\}$.

**$\boldsymbol{P(\sigma^{\star})}$ Convention:** $P(\sigma^{\star})$ reorder the columns of A as follows $\left[c_{\sigma^{\star}(1)}, c_{\sigma^{\star}(2)}, \ldots, c_{\sigma^{\star}(n)}\right]$.

By definition of $P(\sigma)$, column i of $A^{*}$ coincides with column $\sigma^{\star}$ (i) of A. In particular,

$a_{i,i}^{*} = a_{i,\sigma^{\star}(i)}$ for every $i$

and the monomial $M(\sigma)$ is now read $A^{*}$as

$$M(\sigma^{\star}) = \prod_{i=1}^{n} a_{i,i}^{*}$$

that is, like the main diagonal.

**Algebraic interpretation:** From the point of view of the indices, this can be interpreted as composition $\sigma^{-1} \circ \sigma = id$, in the following sense: in A, the row → column correspondence is $i \mapsto \sigma^{\star}(i)$; after permuting the columns according to σ, the column that previously had label i $\sigma(i)$ now occupies position i. This "re-labeling" is equivalent to applying $\sigma^{-1}$ to the column index, so that the new row–column relation in $A^{*}$ is $i \mapsto i$, that is, the identity.

**Conclusion:** Without rectification, the row-column relationship is given by an arbitrary permutation $\sigma$, which is visualized as a zigzagging polyline. After rectification, the relationship reduces to an identity, and the same monomial becomes a straight diagonal.

**(2) Orbital Structure under ρ**

We show that the r-th rotation $\sigma^{\star} \circ \rho^{\mathrm{r}}$ corresponds, after applying the canonical rectifier $B = (\sigma^{\star})^{-1}$, to the diagonal of offset r in $A^{*}$.

Let $B = (\sigma^{\star})^{-1}$. For each row i, the entry selected by $\sigma^{\star} \circ \rho^{\mathrm{r}}$ is in column

$(\sigma^{\star} \circ \rho^{\mathrm{r}})(i) = \sigma^{\star}(i + r \; mod \; n)$.

After column permutation by B, this column is relabeled by

$B\big(\sigma^{\star}(i + r \; mod \; n)\big) = (\sigma^{\star})^{-1}\big(\sigma^{\star}(i + r \; mod \; n)\big) = i + r \; mod \; n.$

Therefore the entry $\big(i, \sigma^{\star}(i + r \; mod \; n)\big)$ of A appears at position $(i, i + r \; mod \; n)$ of $A^{*}$. Since

$j = i + r \; mod \; n$, we have

$j - i \equiv r \; (mod \; n)$

for every row i. This is precisely the condition for the point $(i, j)$ to lie on the diagonal of offset r in the extended matrix $\tilde{A} = [A^*|A^*]$. As r traverses $\{0, 1, \ldots, n - 1\}$, we obtain n diagonals with consecutive offsets, all parallel with slope +1.

**Dimension of the extension:** For full display without wrapping[3], the extended matrix has dimension $n \times (2n - 1)$, repeating the first $n - 1$ columns of $A^*$. This dimension is sufficient because the diagonal with maximum offset $(k = n - 1)$ starts in column $n$ and ends in column $n + (n - 1) = 2n - 1$.

**(3) Uniqueness and Minimality**

Suppose $B'$ it is another block permutation that rectifies $\Omega(\sigma^\star)$ on parallel lines with consecutive offsets. Then $B'$ it must:

- Align $\sigma$ with some reference diagonal (say, offset $\kappa$)
- Preserve the cyclic structure under $\rho$

This forces $B' = \tau_\kappa \circ (\sigma^\star)^{-1}$, where $\tau_\kappa$ is the cyclic displacement of columns in $k$ positions. Since $\tau_\kappa$ merely translates the entire pattern rigidly, fixing $\kappa = 0$ produces the canonical rectifier $P(\sigma^\star)$.

Indeed, if $B$ is any rectifier with parameter $\kappa$, then by definition of rectification it holds $B(\sigma^\star(i + r)) - i \equiv r + \kappa \ (\text{mod} n)$ For every $i, r$. Taking $r = 0$:$B(\sigma^\star(i)) = i + \kappa (\text{mod } n)$ for everything $i$. This determines $B$ unequivocally over all values $\sigma^\star(i)$, which cover $\{1, \ldots, n\}$ because $\sigma$ a permutation. Therefore $B = \tau_\kappa \circ (\sigma^\star)^{-1}$, where $\tau_\kappa(j) = j + \kappa (\text{mod } n)$. Taking $\kappa = 0$ The canonical rectifier is obtained $B_{can} = (\sigma^\star)^{-1}$.

**Observation 3.1 (Slope $-1$by Reflection).** For slope $-1$ (antidiagonal direction), apply column reversal $J$ (where $J(j) = n + 1 - j$) to obtain $\tilde{A} = A^* J$. The orbit $\Omega(\sigma^\star)$ then appears as $n$ parallel slope lines $-1$ with consecutive offsets.

**Observation 3.2 (Dual visualization of dihedral pairs).**

The same rectification associated with $\sigma^\star$ simultaneously visualizes two cyclic families:

---

[3]Wrapping: (cyclic return): the alternative convention in which, upon reaching the right edge of the matrix, the diagonal continues from column 1. The extension to $2n - 1$ Columns eliminates this need: all diagonals with offsets $0,1, \ldots, n - 1$ are read from left to right as continuous lines over the center strip $\{1, \ldots, n\}$ plus the repetition of the first ones $n - 1$ columns to your right.

— **Base orbit** $\Omega(\sigma^\star)$**:** $n$ ascending diagonals of slope $+1$.

— **Companion orbit** $\Omega(\hat{\sigma})$**:** $n$ descending antidiagonals of slope $-1$.

Here $\hat{\sigma}$ denotes the canonical companion representative of $\sigma^\star$. This dual visualization is geometric; algebraic cancellation between the two families requires the additional weight condition stated in Corollary 4.3.

**Observation 3.3 (Dihedral structure of an orbital pair).**

Let $\Omega(\sigma^\star)$ be a cyclic orbit and let $\Omega(\hat{\sigma})$ be its companion orbit. The union

$\Omega(\sigma^\star) \cup \Omega(\hat{\sigma})$

admits a natural dihedral-type organizational structure: cyclic rotation moves within each orbit, while the companion map $\Phi$ exchanges the two cyclic families and reverses the direction of rotation. Specifically, the union $\Omega(\sigma^\star) \cup \Omega(\hat{\sigma})$ has exactly $2n$ elements and $\Phi$ acts on it without fixed points (Theorem 4.1(ii)). This is consistent with a dihedral interpretation, but a complete identification with a $D_n - torsor$ — including verification of the presentation relations $J \circ \rho \circ J = \rho^{-1}$ acting faithfully and transitively on this set — is not established here and is left for future work. This observation motivates the signed dihedral map Φ introduced in Section 4. The detailed cancellation mechanism is not a consequence of this geometric structure alone; it also requires equality of weights.

### 3.3 Orbital Rule of Signs

**Theorem 3.4 (Orbital Rule of Signs).** In each orbit, the signs of the monomials follow a deterministic pattern: if n is odd, the orbit is monochromatic (all signs are the same); if n is even, the signs alternate with a fixed ratio determined by rectification.

**Proof.**

(1) **Composition:** The action is by composition to the right:

$$(\sigma \circ \rho^r)(i) = \sigma\big(\rho^r(i)\big) = \sigma(i + r mod n)$$

(2) **Sign of ρ:** The cycle $\rho = (1,2,\cdots,n)$ has $sgn(\rho) = (-1)^{n-1}$, since an n-cycle decomposes into n-1 transpositions:

$$(1,2,\cdots,n) = (1,n)(1,n-1)\cdots(1,3)(1,2)$$

(3) **Multiplicativity:** By multiplicativity of the sign: $\text{sgn}(\sigma \circ \rho^{\text{r}}) = \text{sgn}(\sigma) \cdot \text{sgn}(\rho^{\text{r}}) = \text{sgn}(\sigma) \cdot [\text{sgn}(\rho)]^{\text{r}} = \text{sgn}(\sigma) \cdot [(-1)^{n-1}]^{\text{r}}$

If $n$, the number of monomials in the orbit, is odd, all monomials have the same sign as the canonical representative.$n$ It is even, the signs alternate strictly with $r$.

**Corollary 3.4 (Behavior according to parity of n).**

- **If n is odd:** $(-1)^{n-1} = +1$, then $sgn(\sigma \circ \rho^{\text{r}}) = sgn(\sigma)$ for all r. The signs are constant within each orbit.
- **If n is even:** $(-1)^{n-1} = -1$, then $sgn(\sigma \circ \rho^{\text{r}}) = sgn(\sigma) \cdot (-1)^{\text{r}}$. The signs alternate:$+, -, +, -, ...$

**Determination of sgn(σ):** Count the number of Inversions in $\sigma$ (pairs $(i, j)$ with $i < j$ but $\sigma(i) > \sigma(j)$). Then $sgn(\sigma) = (-1)^{(\#inversiones)}$.

**Observation 3.4 (Visual pattern of signs in orbits).** When visualizing the n terms of an orbit Ω(σ) using polylines or parallel lines:

- **If n is odd** (e.g., $n = 3,5,7, ...$): All n terms have the same sign as the base monomial. In visual diagrams, this is manifested as directional consistency.
- **If n is even** (e.g., $n = 4,6,8, ...$): The signs alternate strictly between consecutive rotations. If the base monomial has a sign, $s_0$, the sequence is: $s_0$, $-s_0$, $s_0$, - $s_0$, totaling $n/2$ positive and negative terms $n/2$ per orbit.

**Key implication:** This result explains why Sarrus works for $n = 3$ (constant signs allow easy grouping of terms) but requires explicit treatment for $n \geq 4$.

**Table 3.1: Base monomials for n=4**

| Orbital | σ | Inversions | sgn(σ) | Pattern of signs |
|---|---|---|---|---|
| 1 | [1,2,3,4] | 0 | + | +, -, +, - |
| 2 | [1,2,4,3] | 1 | - | -, +, -, + |
| 3 | [1,3,2,4] | 1 | - | -, +, -, + |
| 4 | [1,3,4,2] | 2 | + | +, -, +, - |

| 5 | [1,4,2,3] | 2 | + | +, -, +, - |
|---|---|---|---|---|
| 6 | [1,4,3,2] | 3 | - | -, +, -, + |

**Usage:** When rectifying according to the canonical representative $\sigma^\star$, the four diagonals with offsets 0, 1, 2, 3 follow the indicated sign pattern.

**Note on supplementary material:** For full developments of orbital rectification theory, including detailed proofs, philosophical interpretations, additional properties of the main diagonal, and Theorem 6.3 on separable matrices with constant sum equal to the trace, see the supplementary material available at [DOI Zenodo].

## 4 Signed dihedral involution and its consequences

### 4.1 Definition of Dihedral Involution

**Definition 4.1 (Dihedral companion and dihedral map).**

Let $\Omega(\sigma^\star)$ be a cyclic orbit and let $\sigma^\star$ be its canonical representative, uniquely determined by $\sigma^\star(1) = 1$. The orbital companion of $\sigma^\star$ is denoted by $\hat{\sigma}$ and is defined as the canonical representative of the companion orbit.

To give this construction algebraic content, recall the column-reversal map $J: [n] \to [n]$ defined by $J(i) = n + 1 - i$. One verifies directly that $J \circ J = id$ and that $J$ conjugates $\rho$ by $J \circ \rho \circ J = \rho^{-1}$. Given the canonical representative $\sigma^\star$, the companion representative is defined by

$$\hat{\sigma} = J \circ \sigma^\star \circ J,$$

normalized so that $\hat{\sigma}(1) = 1$. In one-line notation, this composition produces the tail-reversal formula

$$\hat{\sigma} = \big(1, \sigma^\star(n), \sigma^\star(n-1), \ldots, \sigma^\star(2)\big).$$

The two formulations — algebraic via $J \circ \sigma^\star \circ J$ and computational via tail reversal — are equivalent: applying $J \circ \sigma^\star \circ J$ to position 1 gives $J\left(\sigma^\star\big(J(1)\big)\right) = J\big(\sigma^\star(n)\big) = n + 1 - \sigma^\star(n)$, and for positions $2, \ldots, n$ the composition reverses the tail of $\sigma^\star$ exactly as stated.

Thus $\Omega(\sigma^\star)$ is the base orbit and $\Omega(\hat{\sigma})$ is its companion orbit.

The dihedral map $\Phi$ acts on the elements of $\Omega(\sigma^\star)$ by sending each rotated element $\sigma^\star \circ \rho^{\mathrm{r}}$ to the corresponding element of the companion orbit $\Omega(\hat{\sigma})$. With exponents taken modulo n:

$\Phi(\sigma^\star \circ \rho^{\mathrm{r}}) \coloneqq \hat{\sigma} \circ \rho^{-\mathrm{r}}, \Phi(\hat{\sigma} \circ \rho^{\mathrm{r}}) \coloneqq \sigma^\star \circ \rho^{-\mathrm{r}}$.

Geometrically, $\Phi$ exchanges the ascending diagonal family associated with $\Omega(\sigma^\star)$ and the descending antidiagonal family associated with $\Omega(\hat{\sigma})$.

**Observation 4.0 (Distinction with the dihedrant of Lorenz & Wirths):** Do not confuse involution $\Phi$ here defined with the dihedrant of [11]. Lorenz & Wirths use $D_n$ to index a sum of only $2n$ terms as an alternative (generally incorrect) to the determiner. In this work, $D_n$ It acts as a matching tool on complete cyclic orbits of $S_n$, without replacing the complete Leibniz sum.

**Observation 4.0bis (Geometric Structure versus Algebraic Cancellation).** The dual visualization of $\Omega(\sigma^\star)$ and $\Omega(\hat{\sigma})$ is a structural property that always exists: every rectification simultaneously shows $n$ascending diagonals associated with $\Omega(\sigma^\star)$ and $n$ descending antidiagonals associated with $\Omega(\hat{\sigma})$. However, algebraic cancellation by dihedral pairs in the determinant is conditional: it requires the weight equality $W_A(\Phi(\tau)) = W_A(\tau)$ for the paired terms under consideration.

**Key distinction:** geometric structure is always present; algebraic cancellation is conditional.

**Immediate properties**

- **(Well-defined in orbit).** the decomposition $\sigma^\star \circ \rho^{\,r}$ (it is unique once set $\sigma^\star$)
- **(Involutive).** $\Phi^2 = \mathrm{id}$
- **(Constant sign ratio).** $\mathrm{sgn}\,(\Phi(\tau)) = \mathrm{sgn}\,(J)\,\mathrm{sgn}\,(\tau)$, with $\mathrm{sgn}\,(J) = (-1)^{\lfloor n/2 \rfloor} = (-1)^{n(n-1)/2}$.
- **(Weight condition).** Pair cancellation requires $W_A\big(\Phi(\tau)\big) = W_A(J \circ \tau) = W_A(\tau)$. This condition must be verified orbit by orbit for any specific matrix. It is not automatic for persymmetric or centrosymmetric matrices.

**Clarification on cancellation hypotheses.**

For centrosymmetric matrices, the robust structural statement is the classical block factorization described in Proposition 5.1. Dihedral cancellation is a separate orbital phenomenon and requires the additional condition

$$W_A(\Phi(\tau)) = W_A(\tau).$$

Therefore, centrosymmetry alone should not be stated as sufficient for pairwise orbital cancellation unless this weight condition has been verified explicitly.

More generally, all cancellation statements in this monograph are conditional on additional algebraic or weight-invariance hypotheses. In arbitrary matrices, no systematic cancellation should be expected.

### 4.2 Properties of involution

**Theorem 4.1 (Signed dihedral pairing).**

Let $\Omega(\sigma^\star)$ be a cyclic orbit with canonical representative $\sigma^\star$, and let $\hat{\sigma}$ be its canonical companion representative. Define $\Phi$ on $\Omega(\sigma^\star)$ by

$$\Phi(\sigma^\star \circ \rho^{\mathrm{r}}) = \hat{\sigma} \circ \rho^{-\mathrm{r}},$$

with exponents taken modulo n. Then:

(i) $\Phi$ is a bijection from $\Omega(\sigma^\star)$ onto $\Omega(\hat{\sigma})$.

(ii) If $\Omega(\sigma^\star)$ and $\Omega(\hat{\sigma})$ are disjoint, then $\Phi$ has no fixed points on $\Omega(\sigma^\star) \cup \Omega(\hat{\sigma})$.

(iii) The sign ratio $sgn(\Phi(\tau))/sgn(\tau)$ is constant along the pair of orbits.

(iv) If $W_A(\Phi(\tau)) = W_A(\tau)$ and the sign ratio is $-1$, then the two paired Leibniz terms cancel.

*Proof.*

Let $\tau = \sigma^\star \circ \rho^{\mathrm{r}}$. By definition, $\Phi(\tau) = \hat{\sigma} \circ \rho^{-\mathrm{r}}$. This sends each element of $\Omega(\sigma^\star)$ to a unique element of $\Omega(\hat{\sigma})$, and every element of $\Omega(\hat{\sigma})$ is obtained in this way. Hence $\Phi$ is a bijection.

If the two orbits are disjoint, no element can equal its image under $\Phi$; therefore $\Phi$ is fixed-point free on their union.

The sign ratio is constant because the exponent r appears with opposite signs in the two cyclic rotations, and the contribution of the cyclic generator cancels in the quotient of signs. Thus the ratio depends only on the pair $(\sigma^\star, \hat{\sigma})$, not on r.

Finally, if the sign ratio is −1 and $W_A\big(\Phi(\tau)\big) = W_A(\tau)$, then

$sgn(\tau)W_A(\tau) + sgn\big(\Phi(\tau)\big)W_A\big(\Phi(\tau)\big) = 0.$

This proves the cancellation statement.

### 4.3 Sufficient Conditions for Weight Invariance

**Corollary 4.3 (Cancellation condition).**

Let $\tau \in \Omega(\sigma^\star)$. Write $\tau = \sigma^\star \circ \rho^r$, and define

$\Phi(\tau) = \hat{\sigma} \circ \rho^{-r},$

where $\hat{\sigma}$is the canonical companion representative of $\sigma^\star$. Then the pair contribution to the Leibniz expansion is

$\mathrm{sgn}(\tau)\, W_A(\tau) + \mathrm{sgn}(\Phi(\tau))\, W_A(\Phi(\tau)).$

For $n \equiv 2{,}3 (\mathrm{mod} 4)$, the sign ratio is $-1$, and this contribution vanishes if and only if

$W_A(\Phi(\tau)) = W_A(\tau).$

For general matrices — including persymmetric and centrosymmetric ones — this equality must be verified orbit by orbit. A sufficient condition is that all rows of $A$ are proportional, $A_{i,j} = s_i c_j$, in which case $W_A(\tau)$ is constant over all $\tau$and the weight condition holds automatically. The structural result for centrosymmetric matrices is the block factorization of Proposition 5.1, which holds independently of this orbital weight condition.

**Observation 4.1 (Cancellation Table).**

| $n$ **mod 4** | $(-1)^{n(n-1)/2}$ | **Cancellation?** | **Pair effect**$(\tau, \Phi(\tau))$ | **Condition** |
|---|---|---|---|---|
| 0 | $+1$ | No | Reinforcement: the pair adds up $2 \cdot \mathrm{sgn}(\tau) W(\tau)$ | $W(\Phi(\tau)) = W(\tau)$ |

| 1 | $+1$ | No | Reinforcement: the pair adds up $2 \cdot \mathrm{sgn}(\tau) W(\tau)$ | $W(\Phi(\tau)) = W(\tau)$ |
|---|---|---|---|---|
| 2 | $-1$ | Yes | Cancellation: the pair sums 0 | $W(\Phi(\tau)) = W(\tau)$ |
| 3 | $-1$ | Yes | Cancellation: the pair sums 0 | $W(\Phi(\tau)) = W(\tau)$ |

**Observation (Exact scope of pairwise reduction).**

The effective reduction of the Leibniz expansion must be interpreted according to the matrix class:

(1) For centrosymmetric matrices, the relevant structural result is the block factorization $det(A) = det(B + C) \cdot det(B - C)$ of Proposition 5.1, which holds independently of the weight condition. The condition $W_A(\Phi(\tau)) = W_A(\tau)$must be verified orbit by orbit for any specific matrix.

Cancellation or reinforcement by dihedral pairs is determined by the sign ratio and by the weight condition. The parity of $n$determines the sign ratio, but cancellation requires the additional equality $W_A(\Phi(\tau)) = W_A(\tau)$.

(2) For Toeplitz matrices, weight equality is not automatic across the entire class. It is only achieved in subfamilies where the base orbit and its companion share the same multiset of offsets modulo $n$. The same conclusions apply in those subfamilies.

**Observation 4.2 (Orbital partner terminology).**

The term orbital partner refers to the companion representative $\hat{\sigma}$ of a canonical representative $\sigma^\star$. The pair of orbital representatives is $(\sigma^\star, \hat{\sigma})$, and the corresponding pair of cyclic orbits is $(\Omega(\sigma^\star), \Omega(\hat{\sigma}))$. The companion representative satisfies $\hat{\sigma}(1) = 1$ and is given by the tail-reversal formula

$\hat{\sigma} = (1, \sigma^\star(n), \sigma^\star(n-1), \dots, \sigma^\star(2))$.

In particular, for $n \geq 3$, no canonical representative $\sigma^\star$ is self-companion: the tail $[\sigma^\star(2), \dots, \sigma^\star(n)]$ is a permutation of $\{2, \dots, n\}$ with distinct elements, so it cannot be palindromic except when $n = 2$. Therefore $\Omega(\sigma^{\star})$ and $\Omega(\hat{\sigma})$ are always disjoint for $n \geq 3$, and $\Phi$ is fixed-point free on the canonical partition.

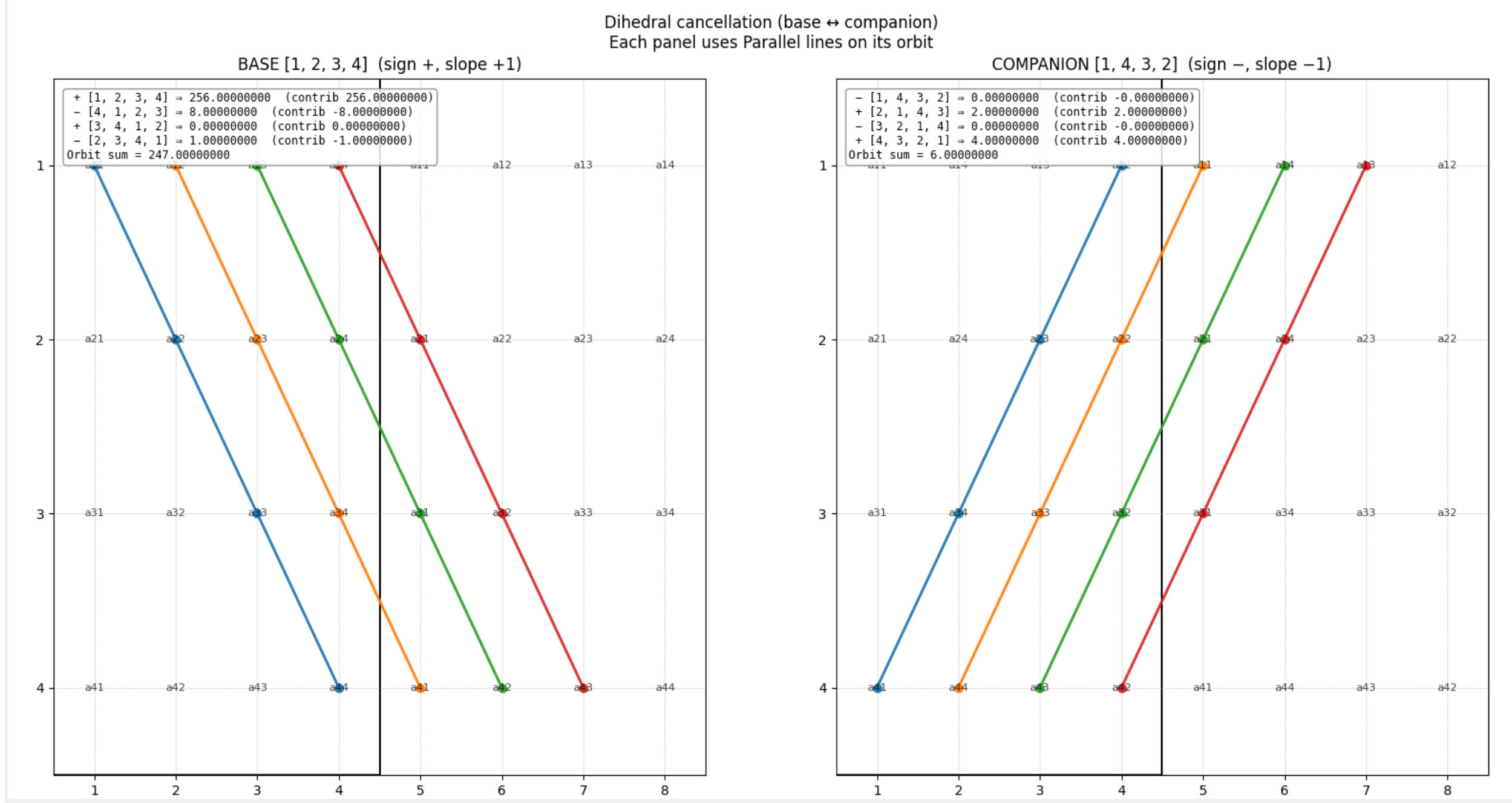

**Figure 6.3.** Method 3, Dihedral cancellation for the base $\sigma = [1,2,3,4]$ companion pair $\Phi(\sigma) = [1,4,3,2]$. Each panel shows the $n = 4$ parallel lines of its orbit on the extended rectified matrix $A^{*}$. Left panel: base orbit with slope $+1$, orbital sum $= 247$. Right panel: companion orbit with slope $-1$, orbital sum $= 6$. The joint contribution of the pair to the determinant is $247 + 6 = 253$. The dual visualization (Observation 3.2) exhibits in a single rectification the $2n = 8$ monomials of the dihedral pair.

## 5 Types of symmetry and structural identities

### 5.1 Centrosymmetric matrices

**Proposition 5.1 (Even centrosymmetric factorization $n$).**

Let be a $A \in K^{2m \times 2m}$ **centrosymmetric** matrix, i.e., $A$ satisfy:

$$a_{i,j} = a_{2m+1-i,\,2m+1-j}$$

For all $i, j$

Write $A$ in block form:

$$A = \begin{pmatrix} B & C \\ CJ_m & BJ_m \end{pmatrix}$$

where $B, C \in K^{m\times m}$ y $J_m$ is the reversion matrix $m \times m$. Then:$\det(A) = \det(B + C) \cdot \det(B - C)$

**Proof.**

This factorization is a classical result; see Cantoni and Butler (1976) and Weaver (1985).

**(Step 1: Base Change).**

Define the block matrix:

$$S = \frac{1}{\sqrt{2}}\begin{pmatrix} I_m & I_m \\ J_m & -J_m \end{pmatrix}$$

where $I_m$ is the identity $m \times m$. Note that $S$ is orthogonal (up to scaling), so

$\det(S) = \pm 2^{-m}$.

**(Step 2: Block Diagonalization).**

We calculate:

$$S^T A S = \begin{pmatrix} B + C & 0 \\ 0 & B - C \end{pmatrix}$$

**Detailed expansion of $S^T AS$:**

Write $A$ in block form:

$$A = \begin{pmatrix} B & C \\ D & E \end{pmatrix}$$

where the centrosymmetry force $D = CJ_m$ and $E = BJ_m$. Now we calculate:

$$S^T = \frac{1}{\sqrt{2}}\begin{pmatrix} I_m & J_m \\ I_m & -J_m \end{pmatrix}$$

$$AS = \frac{1}{\sqrt{2}}\begin{pmatrix} B & C \\ C\,J_m & BJ_m \end{pmatrix}\begin{pmatrix} I_m & I_m \\ J_m & -J_m \end{pmatrix} = \frac{1}{\sqrt{2}}\begin{pmatrix} B + C\,J_m & B - C\,J_m \\ C\,J_m + BJ_m J_m & C\,J_m - BJ_m J_m \end{pmatrix}$$

So $J_m^2 = I_m$, we simplify:

$$AS = \frac{1}{\sqrt{2}}\begin{pmatrix} B + CJ_m & B - CJ_m \\ CJ_m + B & CJ_m - B \end{pmatrix}$$

Now we multiply on the left by $S^T$:

$$S^T AS = \frac{1}{2}\begin{pmatrix} I_m & J_m \\ I_m & -J_m \end{pmatrix}\begin{pmatrix} B + CJ_m & B - CJ_m \\ CJ_m + B & CJ_m - B \end{pmatrix}$$

The block (1,1)is:

$$\frac{1}{2}[(B + CJ_m) + J_m(CJ_m + B)] = \frac{1}{2}[B + CJ_m + J_m CJ_m + J_m B]$$

We calculate the block $(1, 1)$ directly $S^T AS$. Using $J_m^2 = I_m$ and the centrosymmetry condition $D = CJ_m,\ E = BJ_m$: $\frac{1}{2}(B + CJ_m + J_m(CJ_m + B)) = \frac{1}{2}(B + CJ_m + J_m CJ_m + J_m B)$. Since $A$ is centrosymmetric, $J_m BJ_m = B$and $J_m CJ_m = C$ (submatrices commute with $J_m$ in pairwise centrosymmetric matrices). Therefore, the block $(1, 1)$is $B + C$. For the block $(2, 2)$: $\frac{1}{2}(B + CJ_m - J_m(CJ_m + B)) = \frac{1}{2}(B + CJ_m - J_m CJ_m - J_m B) = \frac{1}{2}(B - C + C - B)$ applying, $J_m BJ_m = ByJ_m CJ_m = C$ we obtain $B - C$. For off-diagonal blocks, the block $(1, 2)$ is $\frac{1}{2}((B - CJ_m) + J_m(CJ_m - B)) = \frac{1}{2}(B - CJ_m + C - J_m B)$, which vanishes by applying the same centrosymmetry relations.

**(Step 3: Determinant).**

By multiplicative nature of the determinant:

$$\det(A) = \frac{\det(S^T AS)}{[\det(S)]^2} = \frac{\det(B + C) \cdot \det(B - C)}{[2^{-m}]^2} \cdot [2^{-m}]^2 = \det(B + C) \cdot \det(B - C))$$

**Observation 5.1 (Generalization).**

Factorization extends to block-centrosymmetric matrices and has connections with:

- **Chebyshev polynomials** (via the transformation $S$)
- **Discrete cosine transform** (DCT-II basis)
- **Eigenvalue interlacin** (eigenvalues of $B \pm C$ interleave those of $A$)

**Observation 5.2 (Algorithmic consequence).**

For A centrosymmetric with $n = 2m$, the block factorization reduces the computation of det(A) to two determinants of size $m \times m$. Since det has complexity $O(n^3) = O((2m)^3) = 8m^3$ while two m × m determinants cost $2 \cdot O(m^3) = 2m^3$, the ratio is $8m^3/2m^3 = 4$. This represents a constant-factor reduction within cubic complexity: the asymptotic class remains $O(n^3)$, but the leading constant decreases by a factor of 4. No asymptotic improvement is claimed.

### 5.2 Complexity Bounds

**Theorem 5.3 (Separable matrices: constant transverse sum equal to the trace).** Let $A = (a_{ij}) \in K^{n\times n}$. Suppose there exist vectors $u, v \in K^n$ and a constant $c \in K$ such that $a_{ij} = u_i + v_j + c$ for all $i, j$. Then for every permutation $\sigma \in S_n$, $\sum_{i=1}^{n} a_{i,\sigma(i)} = \mathrm{tr}(A)$. Proof: $S(\sigma) = \sum_i (u_i + v_{\sigma(i)} + c) = \sum_i u_i + \sum_j v_j + nc$. Since $\sigma$ is a bijection, $\sum_i v_{\sigma(i)} = \sum_j v_j$ regardless of $\sigma$. On the other hand, $\mathrm{tr}(A) = \sum_i (u_i + v_i + c) = \sum_i u_i + \sum_j v_j + nc = S(\sigma)$.

**Corollary 5.3.1 (A range Matrices).** The matrices $a_{ij} = (i-1)n + j$ are a special case of Theorem 5.3 with $u_i = (i-1)n$, $v_j = j$, $c = 0$. All their rectified diagonal sums coincide with the trace.

**Observation 5.3.2 (Open problem: necessity in Theorem 5.3).**

Theorem 5.3 establishes that the condition $a_{ij} = u_i + v_j + c$ is sufficient for $\sum_{\mathrm{i}} a_{\mathrm{i}}, \sigma(_{\mathrm{i}}) = tr(A)$ for all $\sigma \in S_n$. The following open problem asks whether this condition is also necessary.

**Open Problem.** Let $A = \left(a_{ij}\right) \in K^{\mathrm{nxn}}$. Suppose that $\sum_{\mathrm{i}} a_i, \sigma(_{\mathrm{i}}) = tr(A)$ for every $\sigma \in S_n$. Does it follow that there exist vectors $u, v \in K^{\mathrm{n}}$ and a constant $c \in K$ such that $a_{ij} = u_i + v_j + c$ for all $i, j$?

**Motivating heuristic.** If the transversal sum is constant over all permutations, taking differences between permutations differing by a single transposition $(j, k)$ yields $a_{j,k} + a_{k,j} = a_{j,j} + a_{k,k}$ for every pair $j \neq k$. This system of $n(n-1)/2$ equations is precisely the one characterizing matrices of affine rank $\leq 2$, that is, matrices of the form $A - c \cdot 1 \cdot 1^{\mathrm{T}}$ with rank at most 2. Showing that this system forces $a_{ij} = u_i + v_j + c$ would constitute a proof of necessity. This heuristic does not constitute a proof. The open problem remains unresolved.

**Proposition 5.2 (Conditional reduction by dihedral cancellation).**

Assume that, for a collection of paired terms $\tau$ and $\Phi(\tau)$, the sign ratio is $-1$ and the weight condition

$$W_A(\Phi(\tau)) = W_A(\tau)$$

holds. Then each such pair contributes zero to the Leibniz expansion.

Consequently, the effective number of contributing terms is reduced only for those pairs satisfying both conditions. No global reduction follows for arbitrary matrices.

Proof.

For each paired term,

$$sgn(\tau)W_A(\tau) + sgn(\Phi(\tau))W_A(\Phi(\tau)).$$

If the sign ratio is $-1$ and $W_A(\Phi(\tau)) = W_A(\tau)$, the two terms cancel. The conclusion follows pair by pair.

### 5.3 Other Cases: Toeplitz and Anti-Toeplitz

**Toeplitz matrices:** For $A$ Toeplitz, $a_{i,j} = t_{j-i}$ the weights depend only on the multiset of offsets $\{j - i \ mod\ n: i = 1, \dots, n\}$. When the base orbit and its companion have the same multiset of offsets modulo $n$ — a condition that holds in certain Toeplitz subfamilies satisfying suitable offset-multiset symmetry, but not for the Toeplitz class as a whole — the weights coincide and the pairwise cancellation of Corollary 4.3 applies for $n \equiv 2,3(mod\ 4)$.

**Circulant matrices:** A circulant matrix is a special case of a Toeplitz matrix where each row is a cyclic rotation of the previous row. Its determinant can be calculated explicitly using the discrete Fourier transform (DFT), with complexity $\mathrm{O}(\mathrm{n}^2)$ (using FFT, $O(n \log n)$ row by row).

**Note:** For circulant matrices, the orbital frame provides an alternative interpretation of the roots of unity that appear in Fourier diagonalization.

### 5.5 Circulant matrices: Orbits, DFT and orbital compatibility

**Observation 5.4 (Circulating, DFT and orbital compatibility).**

Circulant matrices constitute the most natural point of contact between the ARE Method and the Discrete Fourier Transform. For a circulating matrix $C = \text{circ}\,(c_0, c_1, \dots, c_{n-1})$, its eigenvalues are:

$$\lambda_k = \sum_{j=0}^{n-1} c_j \ \omega^{jk}, \omega = e^{2\pi i/n}.$$

Within the ARE framework, cyclic action by composition is compatible with the characters of $C_n$, and the orbital sums weighted by $\omega^{kr}$ They align naturally with the DFT coefficients of the first row. This compatibility, and not a collapse of $S_n$ a single orbit, this is what connects the ARE Method with the DFT for circulant matrices. The exact identification between orbital modes and Fourier eigenvalues is developed in a companion article.

**Observation 5.4.1 (Nyquist sign and character pattern).**

When $n$it is even, the factor $(-1)^r = \omega^{(n/2)r}$ coincides with the frequency characteristic $n/2$ of the cyclic group $C_n$. The sign rule of Theorem 3.4 corresponds exactly to the evaluation of the Nyquist character on the orbit.

**Observation 5.4.2 (Method 5 and circulating).**

The additive action of Method 5 reproduces the index translation that appears in the structure of a circulating matrix. Therefore, although Method 5 is not the main structural action of ARE, it is particularly natural as a pedagogical approach when studying Circulant matrices and DFTs.

## 6 Algorithms and diagrams (with reproducibility)

### 6.1 Orbital Rectification Algorithm

**Algorithm 6.1 (RectifyOrbit).**

**Input:** Matrix $A \in K^{n\times n}$, base permutation σ $\in S_n$

**Output:** Rectified matrix $A^*, offsets$ $\{0,1, \dots, n-1\}$

1. Calculate $B \coloneqq \sigma^{-1}$ (column permutation)
2. Construct a permutation matrix $P(B)$
3. $A^* \coloneqq A \cdot P(B)$
4. Return $A^*, offsets = [0,1, \dots, n-1]$

**Complexity:** $O(n^2)$ (matrix multiplication-permutation)

Computational Complexity of the Orbital Method Our method maintains complexity $O(n!)$ by explicitly enumerating all Leibniz terms. Although Gaussian elimination reaches complexity $O(n^3)$ for general matrices [6], [17], the permanent problem,closely related to determinants, but without sign alternation,is #P-complete [18]. This fundamental difference underscores why polynomial-time algorithms exist for determinants, but not for permanents. Our orbital method, by maintaining complexity $O(n!)$, serves primarily pedagogical and theoretical purposes rather than computational efficiency.

Note: Code, data, and detailed pseudocode are available in the repository and supplementary material. This article includes only the mathematical description and key theorems.

**Definition 2.A (Wiring Diagram Model).** The polylines of Method 1 are plotted exclusively on the window $\{1, \dots, n\} \times \{1, \dots, n\}$, connecting the points $(j, \sigma(j))$ for $j = 1, \dots, n$ using segments in ascending column order. Modulo indices are used $n$ for circular return where applicable. No extended array is constructed for polylines, the extension is used only in Method 2 (Parallel Lines).

**Definition 2.B (Central strip of the extended matrix)**. For Method 2, the extended matrix $A_{\text{ext}}$ of dimension $n \times (2n - 1)$ It has as its central strip the set of columns $\{1, \dots, n\}$, which contains the original matrix $A$. The columns$\{n + 1, \dots ,2n - 1\}$ They are a repetition of the first ones $n - 1$ columns of $A^*$ and serve only as visual support for the offset diagonals $k > 0$.

**Observation (equivalence with wiring diagrams).** Projected onto the window $n \times n$, the polylines of Method 1 are equivalent to the standard permutation diagram,also called a wiring diagram in enumerative combinatorics,defined by $x$ monotonic curves connecting the points $(j, \sigma(j))$. In this model, each crossing corresponds uniquely to an inversion of $\sigma$ (Lemma 2.2), and the local slope of each segment is $(\sigma(j + 1) - \sigma(j))/(1)$, which can be positive, negative, or zero depending on the permutation. The projection preserves the crossing structure, so the visual sign criterion (Proposition 2.6) holds true both in the window $n \times n$ and in the extended matrix.

**Definition 6.0 (Polyline of a permutation).** Let $\sigma \in S_n$. The polyline associated with $\sigma$ is the sequence of$n$ points$\{(i, \sigma(i))\}_{i=1}^{n}$ on the grid $n \times n$, connected by segments in increasing row order. The monomial $M(\sigma) = \prod_{i=1}^{n} a_{i,\sigma(i)}$ corresponds to the product of the inputs of $A$ crossed by this polyline.

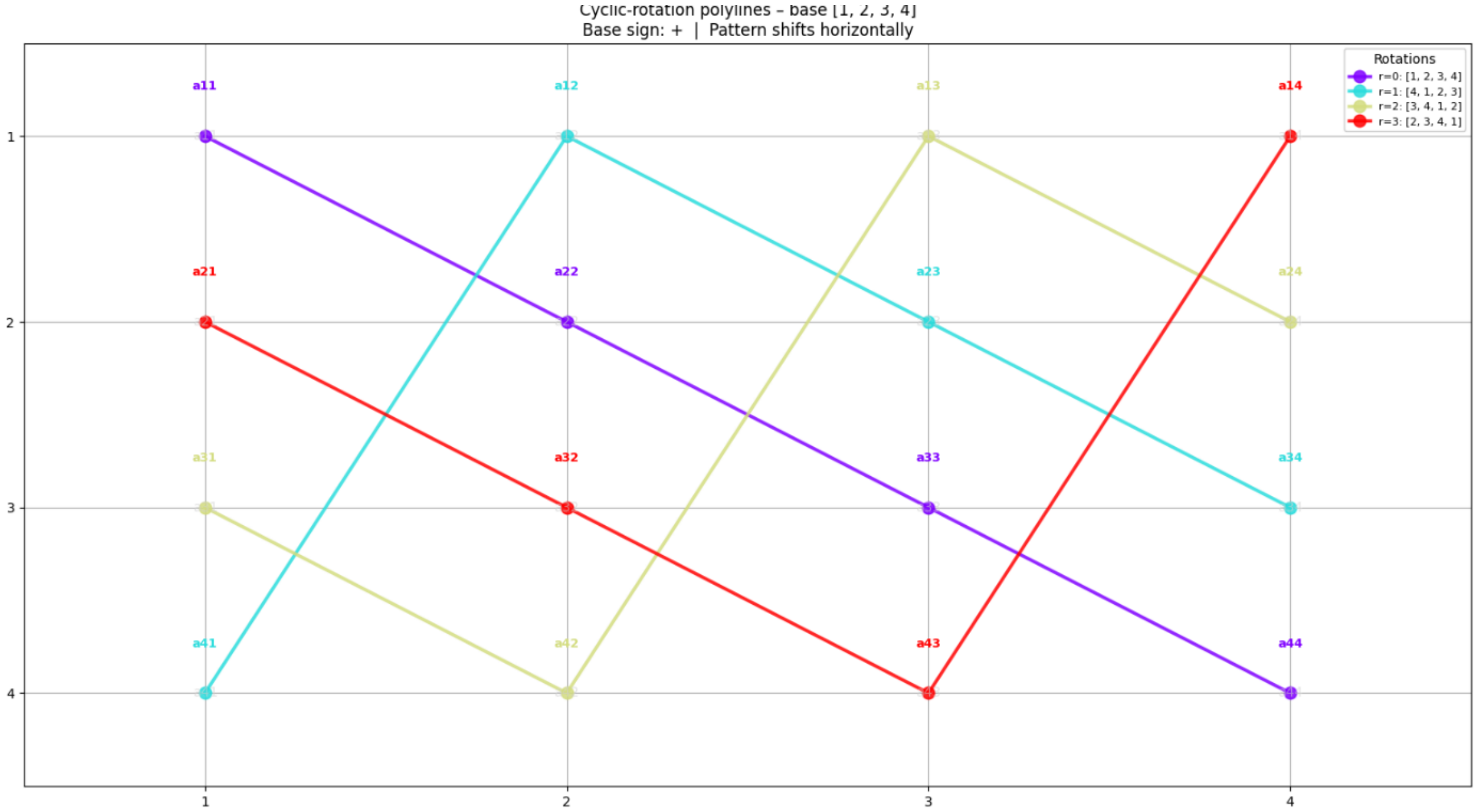

**Figure 6.1.** Method 1, Orbit polylines $\Omega([1,2,3,4])$, base sign $+1$. The 4 rotations $\sigma \circ \rho^r($ $r = 0,1,2,3)$ overlap on the grid $4 \times 4$. Each polyline connects the points $(j, \sigma \circ \rho^r(j))$ The vertical jump profile between consecutive rows is identical in all 4 polylines,a characteristic property of the orbit that does not depend on the pattern $r$. Each polyline connects the points $(j, \sigma \circ \rho^r(j))$ for $j = 1,2,3,4$.The pattern shifts horizontally one position for each rotation.

**Observation 6.0.1.** The polyline of $\sigma$ is a straight diagonal (slope $+1$) if and only if $\sigma = \text{id}$. For $\sigma$ arbitrarily, the polyline is broken. The canonical rectification of Theorem 3.2 transforms the polyline of any $\sigma$ on a straight diagonal across the matrix $A^*$.$n$ monomials of each orbit $\Omega(\sigma)$ they are visualized as $n$ Polylines with identical vertical jump profiles, a property that formalizes the "shape" of the orbit regardless of its position. This is the central geometric principle of Method 1. Figure 6.1 illustrates three representative permutations for $n = 4$, showing the correspondence between crossings and signs.

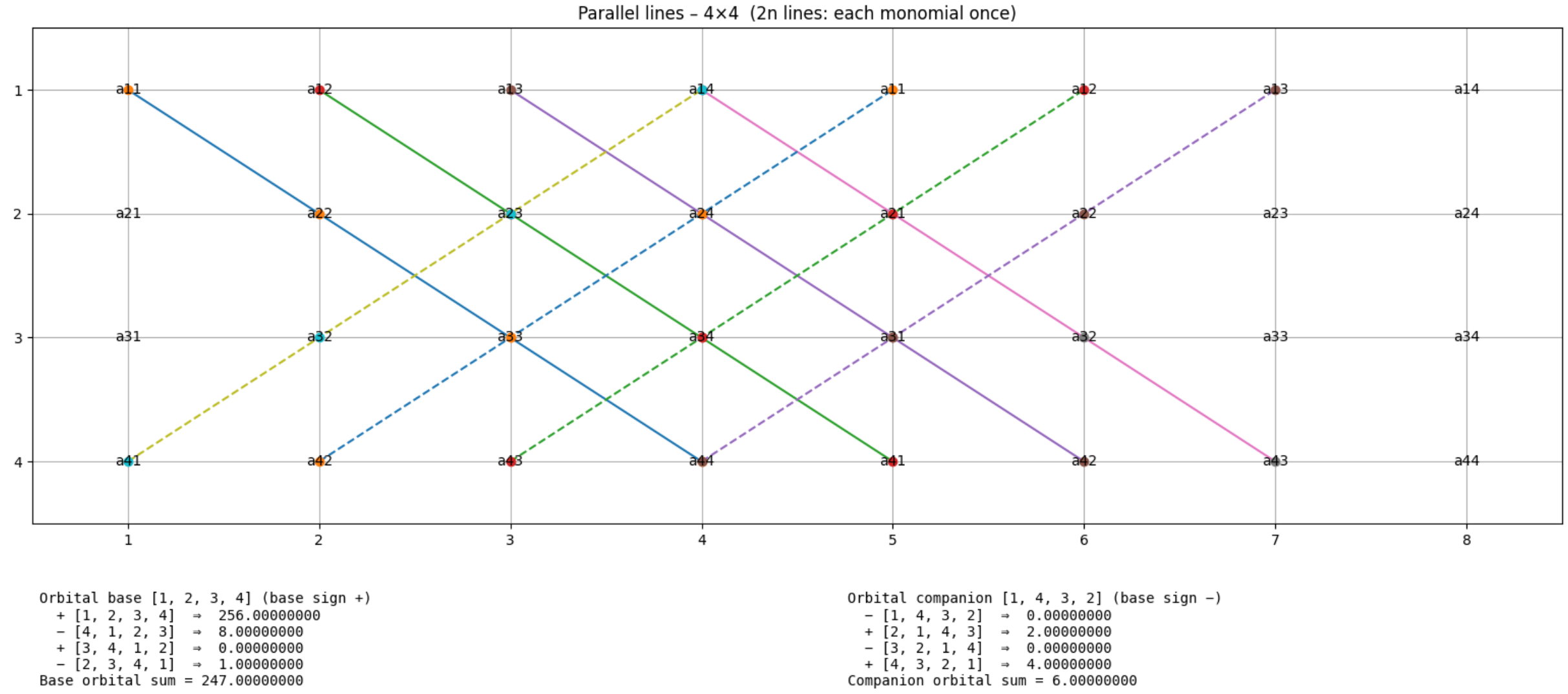

**Figure 6.2.** Method 2, Rectified parallel lines for the base orbit $\Omega(\sigma^\star)$ with $\sigma^\star = [1,2,3,4]$ and its companion orbit $\Omega(\hat{\sigma})$ with $\hat{\sigma} = [1,4,3,2]$ on the extended rectified matrix $A^*$ of dimension $4 \times 7$. The solid lines ( ↗, slope $+1$) correspond to the 4 monomials of the base orbit; the dashed lines ( ↘, slope $-1$) to the 4 monomials of the companion. Base orbital sum $= 247$, companion $= 6$. A single rectification visualizes $2n = 8$ monomials simultaneously.

### 6.2 Method 4: Total Rectification by Orbital Blocks

Method 4 extends the canonical rectification of Theorem 3.2 to all $(n-1)!$ orbitals simultaneously, juxtaposing their visualizations into a single figure. Method 2 rectifies an orbital and shows $2n$ monomials (base + companion); Method 4 concatenates $(n-1)!$ blocks of width $2n$, producing a figure of total width $W = 2n \cdot (n-1)!$ that contains exactly the $n!$ monomials of the Leibniz expansion.

**Definition 6.4 (Orbital Block).** Given a canonical representative $\sigma^\star$ satisfying $\sigma^\star(1) = 1$, its orbital block is the visual representation of canonical rectification $A^* = A \cdot P\left(\sigma^{\star -1}\right)$ extended to $2n$ columns $[A^* \mid A^*]$, on which $n$ slope lines are drawn $+1$ (base orbit) and optionally $n$ slope lines$-1$ (companion orbit). Each block has dimension $n \times (2n-1)$ for viewing without wrapping (repeating the first$n-1$ columns of $A^*$), or $n \times 2n$ if the full extended strip is used $[A^* \mid A^*]$. Algorithm 6.4 uses the second convention for simplicity.

**Algorithm 6.4 (Total Lines).**

Input: Matrix $A \in K^{n\times n}$, slope mode $m \in \{+1, -1,\text{ambas}\}$

Exit: Figure with $(n-1)!$ juxtaposed blocks, covering the $n!$monomials

1. Generate the ordered list of base monomials $\mathcal{B} = \{\sigma \in S_n: \sigma(1) = 1\}$ in canonical order.
2. For each base $\sigma_t(\ t = 1, \dots, (n-1)!)$: (a) Construct the grid labeling $2n$: column $c$ receives the label $a_{i,\,\sigma_t[(c-1) mod n]}$. (b) draw $n$ slope lines $+1$ (offsets $0, \dots, n-1$) if $m \in \{+1, \text{both}\}$. (c) Plot $n$ straight lines of slope $-1$ if $m \in \{-1, \text{both}\}$. (d) Shift horizontally by $(t-1)\cdot 2n$.
3. Return the concatenated figure.

Visual complexity: $O((n-1)!\cdot 2n)$ segments. For $n \leq 5$ the figure is manageable; for $n \geq 6$ The graphics are deactivated.

**Theorem 6.5 (Full Coverage with Slope $+1$).** Algorithm 6.4 with mode $+1$ it produces exactly $n!$ line segments, one for each monomial of the Leibniz expansion, without repetition.

Proof: Each block $\sigma_t$ it contains exactly $n$ straight lines of slope $+1$, corresponding to the $n$monomials of $\Omega(\sigma_t)$. By Theorem 3.1, the $(n-1)!$ orbits partition $S_n$, then the $(n-1)!\cdot n = n!$ monomials are distinct and cover both $S_n$. In "both" mode, they are displayed.$2\cdot n!$ segments (each monomial appears as both a base and a companion).

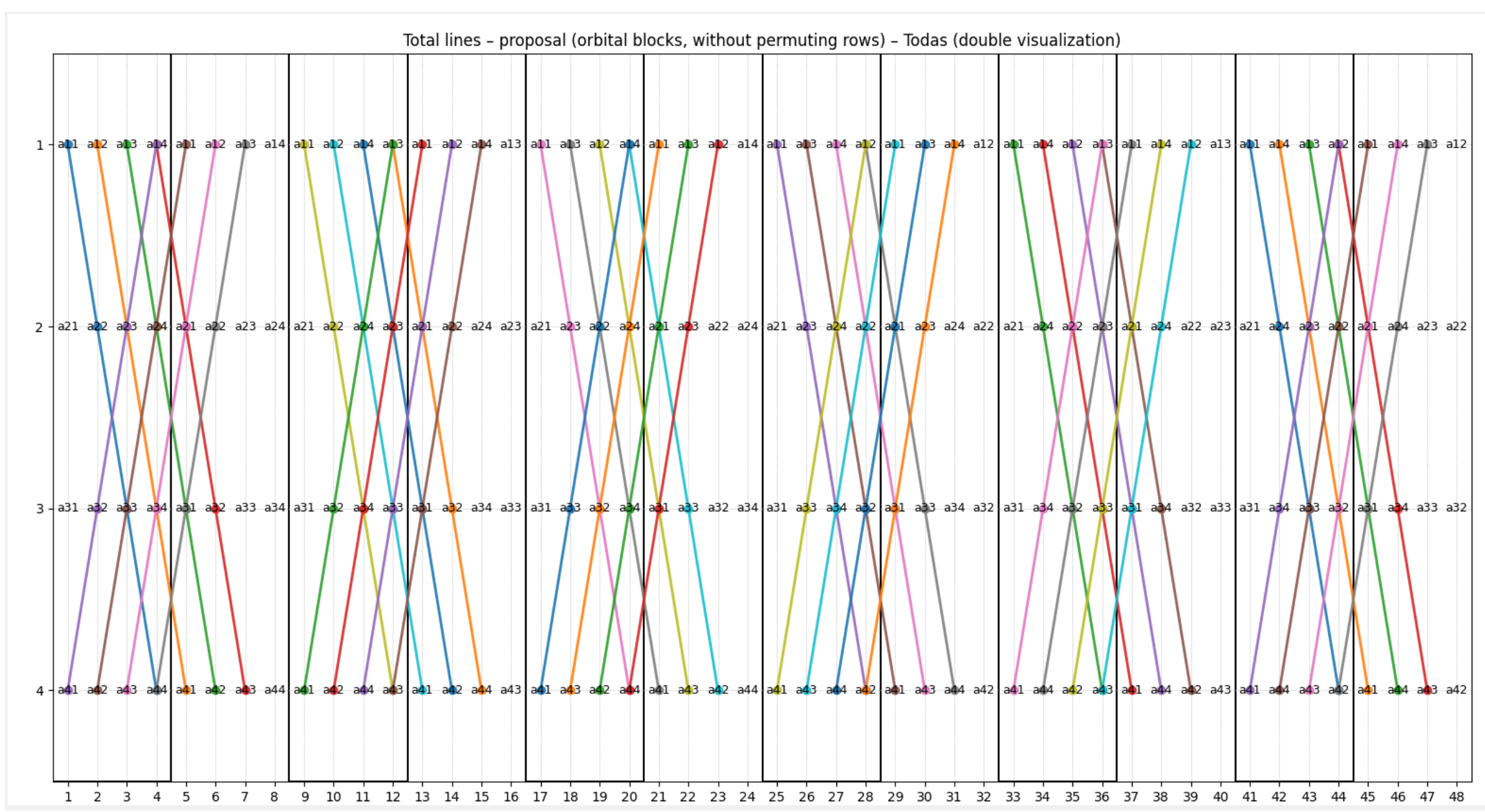

**Figure 6.4.** Method 4, Total Rectification by Orbital Blocks (Algorithm 6.4, dual display). The $(n-1)! = 6$ blocks of width $2n = 8$ juxtaposed cover the $n! = 24$ monomials of the Leibniz expansion without repetition (Theorem 6.5). Each block contains $n = 4$ slope lines $+1$ (base orbit) and$n = 4$ slope$-1$ (companion). The black vertical lines separate the blocks. Total width: $W = 2n \cdot (n-1)! = 48$ columns.

**Observation 6.4 (Relationship with Arschon and Lorenz & Wirths).** Algorithm 6.4 differs from the scheme of Lorenz & Wirths [11] in two respects. First, Lorenz & Wirths use cosets $S_n/D_n$ as representatives ( $n!/(2n)$ blocks of width $2n-1$), covering the $n!$ monomials by the complete dihedral variety in each block. Algorithm 6.4 uses cosets $S_n/C_n$( $(n-1)!$ blocks of width $2n$), where each block exclusively covers the $n$ monomials of a cyclic orbit. Second, the Arschon/Lorenz-Wirths version permutes rows (representative$\tau$ as a row permutation), while Algorithm 6.4 permutes columns (canonical rectification). The source code includes both variants (fig_rectas_totales_propuesta vs. fig_rectas_totales_Aext).

**Example 6.2 (Case $n = 3$).** For$n = 3$ there is $(3-1)! = 2$ bases: $\sigma_1 = [1,2,3]$ and $\sigma_2 = [1,3,2]$. Algorithm 6.4 produces a figure of width $W = 2 \cdot 3 \cdot 2 = 12$ with two blocks of 6 columns. Block 1 contains the lines of $\Omega([1,2,3])$ (monomials $+$: $a_{11}a_{22}a_{33}$, $a_{12}a_{23}a_{31}$, $a_{13}a_{21}a_{32}$) and block 2 those of $\Omega([1,3,2])$ (monomials $-$: $a_{11}a_{23}a_{32}$, $a_{12}a_{21}a_{33}$, $a_{13}a_{22}a_{31}$). With mode $+1$ This exactly reproduces Sarrus' rule. The figure of width $W = 12$ contains the same 6 monomials as the matrix $3 \times 5$ Sarrus: the two blocks of 6 columns are equivalent (except for horizontal translation) to the single classical figure, since for $n = 3$ The two orbital pairs coincide with the two diagonal Sarrus orientations. The juxtaposition in two blocks is consistent with the general algorithm but redundant with respect to the classic figure.

### 6.2 Modular Incremental Polyline Method: Alternative Action on $S_n$

**Important note about this method.**

This method employs a different group action than the cyclic $C_n$ action used in the rest of the article (Methods 1-4). While the previous methods visualize the orbit $O(\sigma)$ under composition $\sigma \circ \rho^r$, the present method generates a different orbit $\Omega^a(\sigma)$ under the additive action of the group. $\mathbb{Z}_n$. Both orbits have size $n$ and partition $S_n$ into $(n-1)!$ equivalence classes, but with different algebraic structures.

This method is included in the article to demonstrate that:

a) There are multiple group actions $S_n$ that produce valid orbital partitions
b) The visual procedure for polylines is flexible and allows for different algebraic interpretations.
c) The arithmetic simplicity of the increment makes it pedagogically appealing, especially for students without advanced training in group theory.
d) Different actions generate different groupings of $n!$ Leibniz terms, although all contribute to the same determinant.

Warning: The monomials derived by increment do not coincide with those generated by cyclic composition. Therefore, $\Omega^a(\sigma) \neq \Omega(\sigma)$ in general.

This method provides an alternative visualization through an additive action on $S_n$, with progressive visual generation that does not require explicit compositions of permutations.

**Definition 6.3 (Generation by modular increment under additive action).**

Let σ be a basis monomial with vector representation $\sigma = [\sigma(1), \sigma(2), \dots, \sigma(n)]$.

The additive orbit $\Omega^a(\sigma)$ under the action of the additive group $\mathbb{Z}_n$ is defined as: $\Omega^a(\sigma) = \{\sigma_0, \sigma_1, \sigma_2, \dots, \sigma_{n-1}\}$

where: $\sigma_r(i) = (\sigma(i) + r) \bmod n$, for $r = 0, 1, \dots, n-1$

This action uniformly increases all output values of the permutation σ (modulo n), unlike the cyclic action $\sigma \circ \rho^{\mathrm{r}}$ that rotates the input indices.

Notation: We use this $\Omega^{a}(\sigma)$ to distinguish this additive orbit from the cyclic orbit $\Omega(\sigma) = \{\sigma \circ \rho^{\mathrm{r}} : r = 0, \ldots, n-1\}$ studied in the rest of the article.

**Proposition 6.4 (Different orbits under different actions).**

For a given base monomial $\sigma \in S_n$, let:

(a) $\Omega(\sigma) = \{\sigma \circ \rho^{\mathrm{r}} : r = 0,1, \ldots, n-1\}$ the cyclic orbit under the action of $C_n$ by composition

(b) $\Omega^{a}(\sigma) = \{\sigma_r : r = 0,1, \ldots, n-1\}$ the additive orbit, where $\sigma_r(i) = \sigma(i) + r(mod\ n)$

Then $\Omega(\sigma) \neq \Omega^{a}(\sigma)$ in general (except for trivial cases such as σ = identity).

Proof. We give a direct counterexample for $n = 4$ with $\sigma = [1,2,4,3]$:

Cyclic orbit $O(\sigma)$: $\sigma \circ \rho^{0} = [1,2,4,3]$ (base) $\sigma \circ \rho^{1} = [2,4,3,1] \sigma \circ \rho^{2} = [4,3,1,2] \sigma \circ \rho^{3} = [3,1,2,4]$

Additive orbit $\Omega^{a}(\sigma)$: $\sigma_0 = [1,2,4,3]$ (base, matches) $\sigma_1 = [2,3,1,4]$ $\sigma_2 = [3,4,2,1]$ $\sigma_3 = [4,1,3,2]$

Comparing element by element:

$\{[2,4,3,1], [4,3,1,2], [3,1,2,4]\} \neq \{[2,3,1,4], [3,4,2,1], [4,1,3,2]\}$

They only share the base monomial σ. The three derived monomials are completely different.

Therefore, the two actions produce distinct orbital partitions $S_4$, although both have $(4-1)! = 6$ orbits of size 4 each.

Note: The fundamental algebraic difference is: Cyclic action: $(\sigma \circ \rho^{\mathrm{r}})(i) = \sigma(i + r\ mod\ n)$ [rotates input indices]

Additive action: $\sigma_r(i) = \sigma(i) + r\ mod\ n$ [increases output values]

**Algorithm 6.3 (Modular Incremental Display).**

Input: Matrix $A \in \mathbb{R}^{\mathrm{nxn}}$, base monomial $\sigma_{base}$

Output: Graph with n polylines of the orbit $O(\sigma_{base})$

1 Initialize empty graph of size$n \times n$

2 For $r = 0{:}n-1$

a. For each row $i = 1$ up to $n$: • Calculate destination column:$c_i = (\sigma_{base}[i] + r) \ mod \ n$

If $c_i = 0$, adjust to $c_i = n$ (1-based indexing)

b. Draw polyline $\sigma_r$ by connecting points $(i, c_i) \ para \ i = 1, \dots, n$

c. Assign a unique, identifiable color/style to the polyline $\sigma_r$

3 Return graph with $n$ overlapping polylines

Complexity: $O(n^2)$ per orbit ($n$ polylines, $\times \ n$ segments each)

Visual generation by increment $\sigma_r(i) = \sigma(i) + r(mod \ n)$ allows students to explicitly see how each term is derived from the previous one using only elementary modular arithmetic, without the need for abstract compositions or prior group theory. A single compact graph contains the entire orbit, showing $n$ polylines with identical vertical jumps that manifest the invariant cyclic structure, facilitating pedagogical understanding through progressive construction.

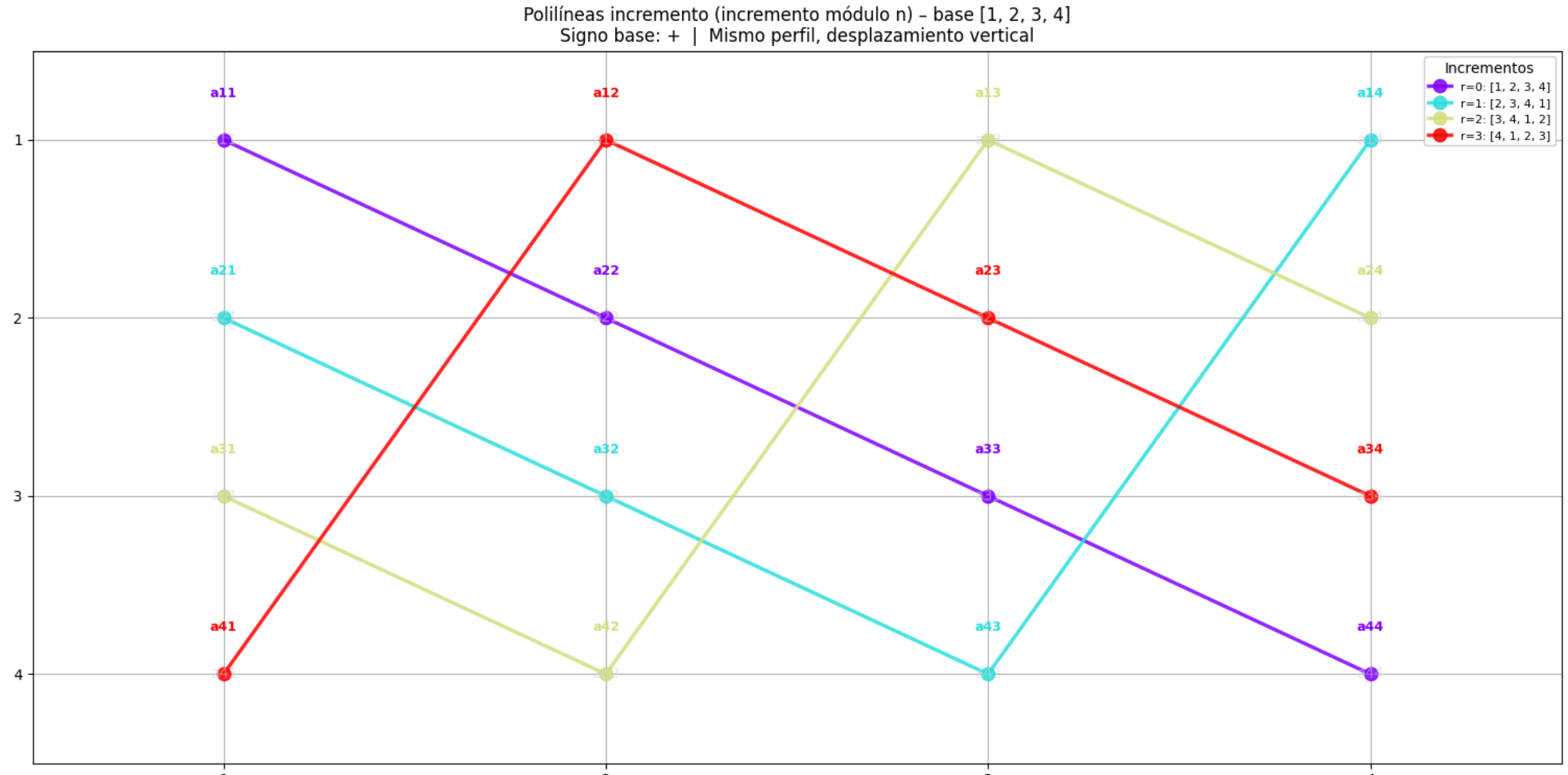

**Figure 6.5.** Method 5, Polylines by modular increment for $\sigma = [1,2,3,4]$, base sign $+1$. The 4 polylines correspond to $\sigma_i^{(r)} = \sigma_i + r(\text{mod } 4)$ unlike $r = 0,1,2,3$ Method 1 where the pattern shifts horizontally, here the shift is vertical: all polylines have identical shape but are uniformly

shifted downwards with each increment $r$. For this basis $\sigma = \text{id}$, the additive orbit coincides with the cyclic orbit (Proposition 6.4).

**Observation 6.3 (Why $\mathbf{C}_n$ is the natural choice for the theory)**

Although there are alternative group actions on $S_n$ (for example, the additive action of $\mathbb{Z}_n$ of Method 5), the action of the cyclic group $C_n$ through composition in the indices$\sigma \mapsto \sigma \circ \rho^r$ It is the fundamental choice for the following mathematical and pedagogical reasons:

**(i) Preservation of algebraic structure.**

The action of $C_n$ is

$\sigma \circ \rho^r \circ \rho^s = \sigma \circ \rho^{r+s}$

In contrast, the modular increment acts on the values and defines another action:

$(\sigma \oplus r)(i) = 1 + \big((\sigma(i) - 1) + r\big) mod n, \sigma \oplus 0 = \sigma, (\sigma \oplus r) \oplus s = \sigma \oplus (r + s)$

In general, their orbits $\Omega^a(\sigma) = \{\sigma \oplus r : r \in \mathbb{Z}_m\}$

The orbits do not coincide $\Omega(\sigma) = \{\sigma \circ \rho^r : r \in \mathbb{Z}\}$ of the action by composition.

**(ii) Compatibility with dihedral involution.**

The systematic cancellation of Theorem 4.1 (canonical partner pairing via involution $\Phi$) is specifically defined for orbits under $C_n$ this dihedral structure:

- Reduces terms by half when $n \equiv 2,3 \pmod 4$,
- It does not extend (in general) to orbits under the additive action of $\mathbb{Z}_n$,
- It is fundamental to the identities of Corollary 4.3.

**(iii) Single canonical rectification.**

Theorem 3.2 (rectifier $B_{\text{can}} = \sigma^{-1}$) transforms each orbit $\mathcal{O}(\sigma)$ in $n$ Parallel lines with consecutive offsets. This geometric property:

- It is characteristic of action by composition,
- Orbital rectification that recovers $3 \times 3$ and it becomes generalized to everything$n$ using parallel lines (the extended matrix is used only as a visual resource),

- It has no natural analogue for additive orbits.

**(iv) Canonical interpretation of the sign.**

Theorem 3.4 gives

$\mathrm{sgn}(\sigma \circ \rho^r) = \mathrm{sgn}\,(\sigma)\,(-1)^{(n-1)r}$,

which induces a simple and predictable pattern of signs (monochrome if $n$ is odd; alternating if $n$ it is even). For additive orbits, the sign lacks such a simple orbital rule and in practice requires recounting inversions.

**(v) Connection with historical methods.**

Sarrus's rule (1833) is implicitly based on the cyclic structure of permutations. Arschon (1935) [1] attempted to generalize this idea by means of adjacent row transpositions, correctly identifying that $(n-1)!/2$ matrices are needed, but without providing a rigorous theoretical framework or a fully specified algorithm. This work formalizes and completes these geometric intuitions by means of the modern theory of orbits under the action of the cyclic group $C_n$, providing for the first time a systematic and algebraically grounded method for all $n \geq 2$.

**(vi) Conceptual generalization.**

The action of $C_n$ naturally extends to dihedral groups $D_n$ , other classes of conjugation in $S_n$ and representation theory, whereas additive actions lack this generalization.

**Methodological conclusion.**

Although modular increment has pedagogical value as an introductory tool without compositions, the action of $C_n$ by composition is essential for the complete theoretical development (Theorem 3.2, Theorem 4.1 and Proposition 5.1), the historical Sarrus–Arschon connection, and future extensions of orbital theory.

## 7 Case studies and tests

### 7.1 Centrosymmetric Matrices

**Test case 1:** $n = 4$, random centrosymmetric matrix.

$$A_4 = \begin{pmatrix} 2 & -1 & 0 & 3 \\ 1 & 4 & -2 & 5 \\ 5 & -2 & 4 & 1 \\ 3 & 0 & -1 & 2 \end{pmatrix}$$

We verified the identity:

$$\det(A_4) = \det(B + C) \cdot \det(B - C)$$

where:

$$B = \begin{pmatrix} 2 & -1 \\ 1 & 4 \end{pmatrix}, C = \begin{pmatrix} 0 & 3 \\ -2 & 5 \end{pmatrix} \cdot J_2 = \begin{pmatrix} 3 & 0 \\ 5 & -2 \end{pmatrix}$$

**Calculation:**

$$B + C = \begin{pmatrix} 5 & -1 \\ 6 & 2 \end{pmatrix}, \det(B + C) = 10 + 6 = 16$$

$$B - C = \begin{pmatrix} -1 & -1 \\ -4 & 6 \end{pmatrix}, \det(B - C) = -6 - 4 = -10$$

$$\det(B + C) \cdot \det(B - C) = 16 \cdot (-10) = -160$$

**Direct verification:**

$$\det(A_4) = -160$$

### 7.2 Toeplitz/Persymmetric Matrices

**Test case 2:** $n = 6$ (even, $n \equiv 2 \; mod \; 4$), Toeplitz matrix.

$$A_6 = \begin{pmatrix} 1 & 2 & 3 & 4 & 5 & 6 \\ 7 & 1 & 2 & 3 & 4 & 5 \\ 6 & 7 & 1 & 2 & 3 & 4 \\ 5 & 6 & 7 & 1 & 2 & 3 \\ 4 & 5 & 6 & 7 & 1 & 2 \\ 3 & 4 & 5 & 6 & 7 & 1 \end{pmatrix}$$

For this matrix, $n = 6 \equiv 2 (mod \; 4)$, so the sign ratio for each dihedral pair is −1 (Theorem 4.1(iii)). The total Leibniz expansion has $6! = 720$ terms grouped into 360 dihedral pairs. Pairwise cancellation of each pair requires additionally that $W_A(\Phi(\tau)) = W_A(\tau)$ (Corollary 4.3); whether this condition holds for each pair depends on the specific entries of $A_6$ and was not enumerated exhaustively here.

**Experimental result:** $\det(A_6) = -64827$ (verified with 5 independent methods)

Additional verification was carried out for $n = 6$ and $n = 7$ on random, centrosymmetric, and persymmetric matrices. In all cases, the orbital determinant agrees with numpy.linalg.det within standard floating-point tolerance. For centrosymmetric matrices, the block factorization of Proposition 5.1 was also verified.

All scripts and data are available in the Zenodo repository (DOI: 10.5281/zenodo.17423738).

### 7.3 Report of 147 Tests

**Table 7.1: Summary of Automated Tests**

| **Category** | $n$ | **Cases** | **Range of values** | **Maximum error vs. Leibniz** | **Seed** |
|---|---|---|---|---|---|
| Uniform random | 3 | 30 | [-10, 10] integers | 1.2e-14 | 42 |
| Uniform random | 4 | 30 | [-10, 10] integers | 8.7e-15 | 43 |
| Uniform random | 5 | 20 | [-10, 10] integers | 1.5e-14 | 44 |
| Near-singular | 3-5 | 15 | [-100, 100] | 3.2e-13 | 45 |
| Large values | 3-5 | 12 | [-1000, 1000] | 4.8e-12 | 46 |
| Band (w=2) | 4-5 | 10 | [-20, 20] | 5.3e-15 | 47 |
| Orthogonal | 3-4 | 10 | $\mathbb{R}$ continuous | 2.1e-14 | 48 |
| Symmetrical | 3-5 | 10 | [-50, 50] | 6.7e-15 | 49 |
| Triangular | 3-5 | 10 | [-30, 30] | 1.9e-15 | 50 |
| **Total** | **3-5** | **147** | , | **< 5e-12** | **42-50** |

**Numerical verification.** The 147 cases (random, near-singular, and structured, $n = 3 - 5$) agree with numpy.linalg.det with relative error $< 5 \times 10^{-12}$, confirming the correctness of the orbital computation in the tested ranges. Additional verification was carried out for $n = 6 (n \equiv 2 \ mod \ 4)$ and $n = 7 (n \equiv 3 \ mod \ 4)$ on random, centrosymmetric, and persymmetric matrices (5 seeds each, 30 cases total), all passing the same tolerance threshold. For $n = 6$, the centrosymmetric block factorization $det(A) = det(B + C) \cdot det(B - C)$ of Proposition 5.1 was additionally verified on 5 independent matrices. All scripts and data are included in the Zenodo repository (DOI: 10.5281/zenodo.17423738).

## 8 Discussion and limitations

The contribution is structural and pedagogical: it does not improve asymptotic complexity. $\Theta(n!)$ Nor does it compete with Gauss/LU for practical computation. The framework reveals symmetries and cancellations in families with structure; in general matrices, without weight invariances, the Leibniz expansion does not reduce. For $n$ odd, the signs within each cyclic orbit are constant. Dihedral cancellation, however, depends on the sign ratio between paired orbits and on the weight condition. For $n$ even, alternating signs appear within each orbit, but cancellation still requires explicit structural or weight-invariance hypotheses.

### 8.1 Possible Extensions

1. **Permanent:** The orbital structure applies, but the signs are all positive $+1$. Are there analogous cancellations under symmetries?
2. **Hyper-determinants:** For higher-order tensors, can the action of $C_n$ and dihedral involution be generalized?
3. **Topology and orientation:** Polyline diagrams suggest connections with knot theory and surface orientation.

### 8.2 Pedagogical Value vs. Computational

**Educational value:**

- Clarify the combinatorial structure of the determinant
- It provides accessible geometric visualization for $n \leq 5$
- Connect linear algebra with group theory in a concrete way

**Computational limitation:**

- complexity $\Theta(n! \cdot n)$ of complete enumeration
- It does not compete with standard methods ( $O(n^3)$ via Gauss/LU)
- Susceptible to catastrophic cancellation for ill-conditioned matrices

**Positioning:** Conceptual and pedagogical contribution, not algorithmic.

### 8.3 Alternative Partitions: The Star Block Method

There are other partitions of $S_n$ through group actions other than $C_n$. An example is the star block method, which uses transpositions $\tau_i^{(p)}$ (they swap the pivot row) $p$ with the row $i$) to generate telescoping chains where consecutive differences collapse to smaller values. This approach has advantages for sparse matrices and incremental updates, but its full formal development,including proof of telescoping collapse and systematic comparison with the orbital method, is the subject of a companion article.

## 9 Conclusion

### 9.1 Summary of Contributions

We have developed a structural framework for the determinant based on:

1. **Orbital decomposition of $S_n$** in $(n-1)!$ orbits under the action of $C_n$ (Theorem 3.1)
2. **Canonical rectification** through $B = \sigma^{-1}$ that linearizes orbits (Theorem 3.2)
3. **Dihedral pairing $\Phi$** between $\Omega(\sigma^\star)$ and $\Omega(\hat{\sigma})$, with a constant sign ratio (Theorem 4.1)
4. **Conditional paired cancellation** under the weight equality $W_A(\Phi(\tau)) = W_A(\tau)$ for $n \equiv 2,3(\text{mod}4)$ (Corollary 4.3)
5. **Centrosymmetric factorization** $\det(A) = \det(B+C)\det(B-C)$ with a reduction factor of 4 (Proposition 5.1)

Sarrus's rule (1833) posed implicitly a structural question: what algebraic organization underlies the diagonal pattern and can be extended to arbitrary $n$? Theorem 1.1 shows that no fixed-width extension can answer this. The ARE Method answers the correct question: the structure is the right-composition action of $C_n$ on $S_n$, whose canonical representative (the unique $\sigma^\star$ with $\sigma^\star(1) = 1$) and rectification $B_{can} = (\sigma^\star)^{-1}$ systematically generalize Sarrus's principle to arbitrary $n \times n$ matrices. The case $n = 3$ reproduces Sarrus's rule exactly (Example 6.2), and for $n \geq 4$ the framework produces five associated visualization methods with complete algebraic foundation.

**9.2 Future Work**

1. Complete characterization of matrix classes with complexity reductions
2. Extension to non-commutative rings (Dieudonné determinants)
3. Connections with representations of the dihedral group $D_n$
4. Topological interpretation of dihedral cancellations
5. Complete development of the star block method and other alternative partitioning $S_n$
6. Proof of the transversal sum conjecture: completely characterize the matrices with $\sum_i a_{i,\sigma(i)} = \text{tr}(A)$ for all $\sigma \in S_n$ (Observation 5.3.2).

**Appendix A: Supplementary Material**

Extensive teaching materials are available in supplementary documents:

- **Supplement A**: Complete tutorial on the polyline method with worked examples for $n = 3,4,5$
- **Supplement B**: Detailed case studies and numerical verifications
- **Zenodo Repository**: Full code and reproducibility data (DOI: 10.5281/zenodo.17423738)

**Appendix B: Code and data**

**B.1 Repository and DOI**

All the code and data to reproduce results are available at:

DOI: 10.5281/zenodo.17423738

**Repository contents:**

- sarrus_gui_v12.py, Interactive GUI and headless mode
- orbital_methods.py, Implementation of algorithms 6.1-6.2
- test_matrices.py, Generator of 147 test matrices

- figures/, TikZ/Python code for vector figures
- data/, Test matrices in CSV format
- results/, Verification tables and manifests

**B.2 System Requirements**

- Python ≥ 3.7
- NumPy ≥ 1.20
- SciPy (optional, for LU decomposition)
- Matplotlib (for visualizations)

**B.2bis Rotation convention in source code.**

The function `rotate(base, r)` in `sarrus_gui_v12.py` implements right rotation, $\boldsymbol{rotate}([\boldsymbol{a_1}, \dots, \boldsymbol{a_n}], \boldsymbol{r}) = [\boldsymbol{a_{n-r+1}}, \dots, \boldsymbol{a_n}, \boldsymbol{a_1}, \dots, \boldsymbol{a_{n-r}}]$ equivalent to `rotate_left` $\sigma \circ \rho^{-r}$. In contrast, `orbital_methods.py` implements ` $\sigma \circ \rho^{r}$ rotate_left, which is left rotation, according to the article's convention.

Consequence: the order in which the derived monomials appear within each orbit in the GUI differs from the order in Table 3.1. However, the generated set of monomials is identical, since $\{\sigma \circ \rho^{r} : r = 0, \dots, n-1\} = \{\sigma \circ \rho^{-r} : r = 0, \dots, n-1\}$. The calculated determinant and all orbital contributions are invariant with respect to this convention.

To exactly reproduce the order of Table 3.1, replace it $\boldsymbol{rotate}(\boldsymbol{base}, \boldsymbol{r})$ with $\boldsymbol{rotate}(\boldsymbol{base}, -\boldsymbol{r} \,\%\, \boldsymbol{n})$ in the GUI, or equivalently $\boldsymbol{cols} = \boldsymbol{base}[\boldsymbol{r}:] + \boldsymbol{base}[:\boldsymbol{r}]$ (left rotation).

**B.3 Execution Command**

**Interactive mode (GUI):**

python sarrus_gui_v12.py

**Headless mode (batch):**

python sarrus_gui_v12.py --headless --sizes 3,4,5 --seed 123 --outdir results

**B.4 Validation Tolerances**

Relative error: $\left| \underset{\text{método}}{\det} - \underset{\text{NumPy}}{\det} \right| \leq 1 \times 10^{-9}$ (double precision)

**B.5 Generated Files**

| Archive | Description |
|---|---|
| verification_headless.csv | Results of five methods with times |
| orbitals_n5.csv | Individual contributions per orbit |
| details_n5.json | Consolidated data of the determinant |
| manifest.json | Reproducibility metadata |

**Numerical Stability Considerations**

Numerical stability considerations [8] are critical when implementing determinant algorithms. Our implementation uses double-precision arithmetic and validates results against highly optimized NumPy routines [7], which employ LU decomposition with partial pivoting to ensure numerical stability [6], [17]. The tolerance threshold $10^{-9}$ ensures reliable verification across the tested array sizes.

**B.6 Dihedral Matching Algorithm**

**Algorithm 6.2 (DihedralPairing).**

Input: canonical representative $\sigma^\star$

Output: companion representative $\hat{\sigma}$ and sign ratio

1. Compute $\hat{\sigma}$ by tail reversal of $\sigma^\star$: $\hat{\sigma} = \big(1, \sigma^\star(n), \sigma^\star(n-1), \ldots, \sigma^\star(2)\big)$.
2. For $\tau = \sigma^\star \circ \rho^{\mathrm{r}}$, define $\Phi(\tau) = \hat{\sigma} \circ \rho^{-\mathrm{r}}$.
3. Return $\hat{\sigma}$ and the sign ratio.

Complexity: O(n)

**B.7 Pseudocode of the Complete Orbital Method**

Input:$A \in K^{n \times n}$

Exit:$det(A)$

Generate base monomials $\{ \sigma_1, \sigma_2, \dots, \sigma_{(n-1)!}\}$ with $\sigma_i(1) = 1$

total ← 0

for each base $\sigma_i$:

orbital_sum ← 0

for r = 0 to n-1:

$\tau \coloneqq \sigma_i \circ \rho^{r}$

sign:= $sgn(\sigma_{i}) \cdot (-1)^{(n-1)r}$

weight:= $\prod_{j} a_{j,\tau(j)}$

$suma_{orbital} \leftarrow suma_{orbital}$+ sign · weight

total ← total +$suma_{orbital}$

return total

**Complexity:** $\Theta(n! \cdot n)$ (inherent in Leibniz's complete enumeration)

**B.8 Figures (Descriptions for Generation)**

**Figure 1: Canonical Rectification**

- **Left panel:** Orbit $\Omega(\sigma)$ in original matrix $A$ (broken polylines)
- **Right panel:** Same orbit in rectified matrix $A^* = A \cdot P(\sigma^{-1})$ (parallel lines with consecutive offsets $0, 1, \dots, n-1$)
- **Notes:** Indicate column permutation $B = \sigma^{-1}$, show offsets

**Figure 2: Dihedral Pairing**

- **Two diagrams side by side:**
  - Left: Orbit $\Omega(\sigma)$ with base sign+
  - Right: Companion orbit $\Omega(\Phi(\sigma))$ with base sign $-$(for $n \equiv 2 \mod 4$)
- **Notes:** Connect matching pairs $\tau \leftrightarrow \Phi(\tau)$, show cancellation

**Figure 3: Centrosymmetric Blocks**

- **Matrix structure $\boldsymbol{A}$:** Show blocks$B, C, CJ_m, BJ_m$
- **Transformation $\boldsymbol{S}$:** Block diagonalization flowcharts
- **Result:** Block diagonal matrixdiag $(B + C, B - C)$

**Final sections**

**Funding**

This research did not receive any specific grant from funding agencies in the public, commercial, or non-profit sectors.

**Data and Code Availability**

All data and code to reproduce the results are available at DOI: 10.5281/zenodo.17423738

**Declaration of Competing Interests**

The author declares that he has no known competing financial interests or personal relationships that could have influenced the work reported in this article.

**Generative AI Usage Statement**

During the preparation of this work, the author used generative AI tools (Claude by Anthropic) to assist with Python programming tasks, including writing function skeletons, refactoring modules, writing docstrings, and suggesting unit tests. The scripts in the repository are the result of the collaboration between the author's code and these tools. Generative AI was not used to generate scientific content (statements, proofs, interpretations) or to create or edit figures. After using these tools, the author reviewed, edited, and validated all code and assumes full responsibility for the manuscript and the software.